\documentclass[12pt]{article}

\usepackage{epsfig,amsmath,amsfonts,amssymb,latexsym,color,amsthm,hhline,multirow,subfigure,graphicx}

\newcommand{\E}{{\mathbb E}}

\newcommand {\PP}{{\mathbb P}}
\newcommand{\sss}{\scriptscriptstyle}

\begin{document}

\title{Birds of a feather or opposites attract - effects in network modelling}\parskip=5pt plus1pt minus1pt \parindent=0pt
\author{Maria Deijfen\thanks{Department of Mathematics, Stockholm University; {\tt mia@math.su.se}} \and Robert Fitzner\thanks{Department of Mathematics and Computer Science, Eindhoven University of Technology; {\tt math@fitzner.nl}}}
\date{July 2017}
\maketitle

\begin{abstract}
\noindent We study properties of some standard network models when the population is split into two types and the connection pattern between the types is varied. The studied models are generalizations of the Erd\H{o}s-R\'{e}nyi graph, the configuration model and a preferential attachment graph. For the Erd\H{o}s-R\'{e}nyi graph and the configuration model, the focus is on the component structure. We derive expressions for the critical parameter, indicating when there is a giant component in the graph, and study the size of the largest component by aid of simulations. When the expected degrees in the graph are fixed and the connections are shifted so that more edges connect vertices of different types, we find that the critical parameter decreases. The size of the largest component in the supercritical regime can be both increasing and decreasing as the connections change, depending on the combination of types. For the preferential attachment model, we analyze the degree distributions of the two types and derive explicit expressions for the degree exponents. The exponents are confirmed by simulations that also illustrate other properties of the degree structure.

\vspace{0.3cm}


\vspace{0.2cm}

\end{abstract}
\section{Introduction}

The fact that people that are similar in some respect (sex, opinion, class etc) often tend to group together is a well-known phenomenon in sociology - ''birds of a feather flock together''. It is referred to as \emph{homophily}. In other contexts it might, on the contrary, be the case that people with different features are drawn to each other - ''opposites attract''. This is referred to as \emph{heterophily}. See \cite{Bof, Rogers_Bhowmik} and the references therein for details. The last decade has seen the growth of many online social networks and the same phenomena may of course occur there; see e.g.\ \cite{dating,online} for studies on an online dating community and Facebook, respectively. Homophily/heterophily is also present in many other types of networks, for instance the WWW, where the types of webpages may affect the link structure, and biological networks, where specific properties of a gene or protein may be important for interactions.

The aim of the present paper is not to introduce fundamentally new models, but to investigate how well-known standard models behave when the particular phenomenon of homophily/heterophily is incorporated and tuned. We restrict to the very simple case with two types, that is, each member of the population belongs to precisely one of two groups, e.g.\ men/women, left/right-wing in the political spectrum, higher education/no higher education, business webpage/personal webpage, company/private email domain etc. All models that we consider can easily be extended to incorporate several types, but the analysis becomes somewhat more elaborate.

We study two-type versions of the Erd\H{o}s-R\'{e}nyi graph \cite{RG_Boll, RG_Svante}, the configuration model \cite{Boll_conf,MR-95, MR-98} and a simple preferential attachment model \cite{BA,BRST}. For the Erd\H{o}s-R\'{e}nyi graph and the configuration model, the interest revolves around the threshold for the occurrence of a giant component and the size of such a component above the threshold. Preferential attachment models typically do not exhibit phase transitions in the component structure -- indeed, our version gives a connected graph for all parameter values -- but the focus there is instead on various properties of the degree distribution. Computer simulations are used to illustrate results and investigate the models further. The source code is available online; see \cite{code}.

\noindent \textbf{The Erd\H{o}s-R\'{e}nyi graph.} In the Erd\H{o}s-R\'{e}nyi graph, each vertex pair in a set of $n$ vertices is connected independently with probability $q$. When $q=\alpha/n$ for some $\alpha>0$, the model exhibits a phase transition at $\alpha=1$: if $\alpha>1$ there is a unique giant component of order $n$, while if $\alpha\leq 1$ the largest component is of strictly smaller order (logarithmic in $n$ for $\alpha<1$). In the two-type version, each vertex is independently classified as type 1 or type 2 with probability $p_1$ and $p_2=1-p_1$, respectively. The edge probability between two type $i$ vertices ($i=1,2$) is $\min\{\alpha_i/n,1\}$ and the edge probability between two vertices of different type is $\min\{\beta/n,1\}$, where we assume that $\beta>0$. Each possible edge in the graph is included independently. This model is a simple special case of the very general model in \cite{BJR} and was first studied in \cite{Sod}. The model setup here is hence not new, but our contribution lies in analyzing the properties of the above specific instance of the model with respect to the parameters $\alpha_1$ and $\alpha_2$ (controlling homophily) and $\beta$ (controlling heterophily).

In Section 2, we give an expression for the threshold for the occurrence of a giant component and investigate how it varies with the homophily/heterophily parameters. For our examples, we typically keep the proportions of the types and their expected degrees fixed. This means, in particular, that the edge density in the graph is also fixed. We then vary the proportion of edges connecting vertices of the same and different type, respectively. Generally, we see that the threshold parameter decreases as vertices become more prone to connect to vertices of opposite type. This is because the type with the smaller expected degree then gets a larger influence of the graph.

The size of the largest component is studied by aid of simulations and shows different behavior depending on the combination of types. If a subcritical type (that is, a type that is subcritical on its own) and a supercritical type are combined, and the combination is such that the critical parameter eventually drops below 1, then the component size decreases as the mixing between the types increases. If the subcritical type is not able to drive the graph to subcriticality -- that is, if the critical parameter stays above 1 -- then the component size may, on the contrary, increase as edges shift to connecting vertices of different type. This is because superfluous edges of the supercritical type may then be used to include vertices of the subcritical type in the giant component without loosing vertices of the supercritical type; see Section 2. Combining two supercritical types, there is a sharp increase in the component size immediately when the two types start to mix, caused by two giant components merging, but the component size is then not affected much as the connection pattern changes further.

\noindent \textbf{The configuration model.} The Erd\H{o}s-R\'{e}nyi graph and its above generalization give rise to Poissonian degree distributions. Many empirical networks however exhibit more heavy-tailed degree distributions. The configuration model \cite{Boll_conf,MR-95, MR-98} provides a way of generating a graph with an arbitrary degree distribution. Here we describe a two-type version of it, with potentially different degree distributions for the two types and where the fractions of edges between different combinations of types can be tuned. A similar model is used in \cite{Ball_Sirl} in the context of multi-type epidemics.

As before, the number $n$ of vertices is fixed and each vertex is assigned type 1 or 2 with probability $p_1$ and $p_2$, respectively. Let $F_1$ and $F_2$ be two probability distributions with support on the non-negative integers and finite means. For $i=1,2$, independently equip each type $i$ vertex with a random number of half-edges according to the distribution $F_i$. Each half-edge at a type 1 vertex is independently labeled 1 or 1' with probability $\xi_1$ and $1-\xi_1$, respectively, and similarly each half edge at a type 2 vertex is independently labeled 2 or 2' with probability $\xi_2$ and $1-\xi_2$, respectively. Half-edges with label 1 are then paired to each other uniformly at random, and similarly for half-edges of type 2. This creates edges that connect vertices of the same type. Half-edges with label 1' are randomly paired to half-edges with label 2', creating edges connecting vertices of different type. Half-edges that remain after the pairing procedure are erased.

In Section 3, we give a condition that guarantees that the degree distributions are asymptotically not affected by the erasing of un-paired half-edges, and we derive an expression for a threshold parameter for the occurrence of a giant component. Similarly to the analysis for the Erd\H{o}s-R\'{e}nyi graph, we then investigate how the threshold changes as the allocation of the edges between/within the vertex types varies. The conclusion is analogous: a larger element of heterophily decreases the threshold parameter. Also the conclusions for the size of the largest component are similar to the Erd\H{o}s-R\'{e}nyi case: for combinations of a subcritical and a supercritical type the size can be either increasing or decreasing depending on the specific combination, while combining two supercritical types leads to a roughly constant component size. There are also cases where the change in the component size is not monotone as the connection pattern changes; see Section 3.

\noindent \textbf{Preferential attachment.} Finally, we consider a version of the basic preferential attachment model \cite{BA,BRST}, with two vertex types, and where a new vertex chooses which existing vertex to connect to based on both type and degree. At time $t=1$, the graph consists of two vertices, one of type 1 and one of type 2, connected by an edge. At each integer time $t\geq 2$, a new vertex with one edge attached to it arrives in the graph. With probability $p_1$ the new vertex is type 1 and with probability $p_2=1-p_1$ it is type 2. A type $i$ vertex that arrives in the graph connects to another type $i$ vertex with probability $\theta_i$ and to a vertex of opposite type with probability $1-\theta_i$. Finally, when a new vertex has decided which type to connect to, it chooses a vertex of that type proportionally to degree; see Section 4 for further details.

As mentioned, this model gives rise to a connected graph, and questions about component size are therefore not relevant. Instead, we focus on the degree distribution. In Section 4 we derive an expression for the expected asymptotic degree distributions for the types and conclude that they both obey power-laws. Specifically, the fraction of type $i$ vertices with (total) degree $k$ obeys a power-law with exponent
$$
\tau_i:=2+\frac{p_i}{p_i\theta_i+p_{i^c}(1-\theta_{i^c})},
$$
where $i^c$ denotes the complementary type. It follows e.g.\ that, in a population where vertices connect only to vertices of opposite type, the exponent for type $i$ is $2+\frac{p_i}{p_{i^c}}$, so that the type that occupies the smallest fraction of the population has a heavier degree tail. When both types connect only to type $i$, we get $\tau_i=2+p_i$ and $\tau_{i^c}=\infty$, so that type $i$ has a heavier tail than in standard preferential attachment. In Section 4, we analyze the degree distributions further for some specific instances of the model and compare our results to simulations. We also investigate the number of neighbors of a specific type and the dependence between these type-specific degrees.

\noindent \textbf{Previous work.} Although the phenomena of heterophily/homophily are well-known from empirical work (see e.g.\ the references above) there has been relatively little work aimed at including them in probabilistic network models. Graph models where weights are associated with the vertices and the edge probabilities determined by the weights are frequent in the probability literature; see e.g. \cite{GenGraph, BJR, DEHH,NR}. The weights could of course be interpreted as types, but typically the edge probabilities are increasing functions of the weights and do not take the difference between weights into account. In such a setting the weights represent fitnesses of the vertices rather than types without numerical importance. The very general model in \cite{BJR} however is well suited for modeling homophily/heterophily and, as mentioned, the two-type version of the Erd\H{o}s-R\'{e}nyi graph above is a special case of this model. The specific analysis in Section 2 of the effect of tuning the connection pattern has, to our best knowledge, not appeared before, and we have therefore included it here.

The model in \cite{attribute} is specifically aimed at constructing graphs where the edge probabilities are determined by attributes associated with the vertices. However, it deals mainly with the situation when the number of attributes is very large and, due to restrictions on the parameter space, it also does not allow for altering the connection pattern in the way that we are interested in. The model in \cite{attribute} has close links to previous work in the statistical literature on social networks aimed at including various type of information about individuals in network models; see e.g.\ \cite{latent, block}. However, the focus there is typically on estimating parameters rather than deriving asymptotic mathematical properties. In this context, we also mention the planted partition model, where the population is partitioned in two classes with varying link propensity between and within the classes, and the goal is to reconstruct the underlying partition of the vertex set given the outcome of the graph; see e.g.\ \cite{planted} and references therein.

The rest of the paper is organized so that the above models are analyzed in more detail in Section 2-4, and Section 5 then contains conclusions and some suggestions for further work.

\section{The Erd\H{o}s-R\'{e}nyi graph}

We first deal with the two-type Erd\H{o}s-R\'{e}nyi graph described above. Recall that a vertex is type $i$ ($i=1,2$) with probability $p_i$, the edge probability between two type $i$ vertices is $\min\{\alpha_i/n,1\}$ and the edge probability between two vertices of different types is $\min\{\beta/n,1\}$ for $\beta>0$. We use the notation $i^c$, where $i^c=1$ when $i=2$, and $i^c=2$ when $i=1$.

Write $D(n)$ for the degree of vertex 1 and $D^{\sss (i)}(n)$ for the type $i$ degree of vertex 1, that is, the number of type $i$ neighbors of vertex 1 (note that the vertices are exchangeable so it is enough to consider vertex 1). Furthermore, let $T$ denote the type of vertex 1. It follows from \cite[Theorem 3.13]{BJR} that $D(n)$ converges in distribution to a mixed Poisson variable $D$ as $n\to \infty$, with
$$
\mu_i:=\E[D|T=i]=\alpha_ip_i+\beta p_{i^c}.
$$
The type $i$ degree also converges in distribution to a mixed Poisson variable $D^{\sss (i)}$, with
$$
\E[D^{\sss (i)}|T=i]=\alpha_ip_i\quad \mbox{and}\quad \E[D^{\sss (i)}|T=i^c]=\beta p_{i^c}.
$$
The model exhibits a phase transition in the sense that there exists a critical parameter $\lambda_c=\lambda_c(p_1,p_2,\alpha_1,\alpha_2,\beta)$ such that, if $\lambda_c>1$, there is a unique giant connected component whose size is of order $n$ while, if $\lambda_c\leq 1$, the size of the largest connected component is of strictly smaller order (logarithmic in $n$ for $\lambda_c<1$); see \cite[Theorem 3.1 and 3.12]{BJR}. For $i,k\in\{1,2\}$, let $m_{ik}$ denote the expected number of type $k$ neighbors of a type $i$ individual. The critical parameter $\lambda_c$ is given by the largest eigenvalue of the matrix $M=\{m_{ik}\}$ (the result \cite[Theorem 3.1]{BJR} is expressed in terms of a certain integral operator and covers a more general setup, but the expression is equivalent to the largest eigenvalue in our case). The result is based on theory for multi-type branching processes. We have
$$
M=\begin{pmatrix}
\alpha_1p_1 & \beta p_2 \\
\beta p_1 & \alpha_2 p_2
\end{pmatrix}
$$
and standard calculations give
\begin{equation}\label{eq:ERl_allm}
\lambda_c=\frac{\alpha_1 p_1+\alpha_2 p_2}{2}+\sqrt{\left(\frac{\alpha_1p_1+\alpha_2p_2}{2}\right)^2+ p_1p_2(\beta^2-\alpha_1\alpha_2)}.
\end{equation}
For $\lambda_c>1$, the asymptotic size of the giant component can be expressed in terms of generating functions related to the branching process. The expressions however become rather involved, and we will instead investigate the size by aid of simulations.

How is the critical parameter and the size of the giant component affected by the level of homophily/heterophily in the graph? The homophily level is controlled by the parameters $\alpha_1$ and $\alpha_2$: as $\alpha_i$ increases, type $i$ vertices become more prone to create connections to other type $i$ vertices. Similarly, heterophily is controlled by the parameter $\beta$: increasing $\beta$ makes the two types more inclined to bond with each other. Of course, if one of these parameters increases while everything else is fixed, the number of edges in the graph also increases. To enable a fair comparison of the properties of the graph, we will keep the expected degree in the graph (and thereby also the expected number of edges) fixed, in fact we will even keep the expected degree $\mu_i$ for each type fixed. This means that, while increasing e.g. $\beta$ we will decrease $\alpha_1$ and $\alpha_2$ so that some edges are shifted from connecting vertices of the same type to connecting vertices of different type. Note at this point that the expression (\ref{eq:ERl_allm}) for $\lambda_c$ can be written as a function of $\beta$ and $\mu_i$ ($i=1,2$):
$$
\lambda_c=\frac{1}{2}\left(\mu_1+\mu_2-\beta+\sqrt{(\mu_1-\mu_2)^2+\beta^2+2\beta(\mu_1-\mu_2)(p_1-p_2)}\right).
$$
It is straightforward to confirm that the derivative with respect to $\beta$ of this function is negative (while the second derivative is positive). The function is hence decreasing (and convex) in $\beta$, that is, increasing  $\beta$ decreases the critical parameter (and the decrease is largest for small values of $\beta$). The reason is that, when the types start to mix more, the type with the smaller expected degree gets a larger influence of the graph. We will clearly see this in Examples 2.2-2.4 below.

The size of the largest connected component will be studied by aid of computer simulations. Throughout, graphs of size $n=10,000$ are simulated and the pictures show averages of 100 realizations per parameter value. The code is written in Java and is available online; see \cite{code}.

\noindent \textbf{Example 2.1: The symmetric case.} First consider the totally symmetric case $p_1=1/2$ and $\alpha_1=\alpha_2=\alpha$. The expected degree of both types (as well as the average expected degree) equals $(\alpha+\beta)/2$, and the expression for $\lambda_c$ reduces to $(\alpha+\beta)/2$. The value of the critical parameter is hence not affected if the expected degree(s) are fixed.\hfill$\Box$

\begin{figure}[p!]
\centering \mbox{\subfigure[The critical parameter.]{\includegraphics[height=5cm]{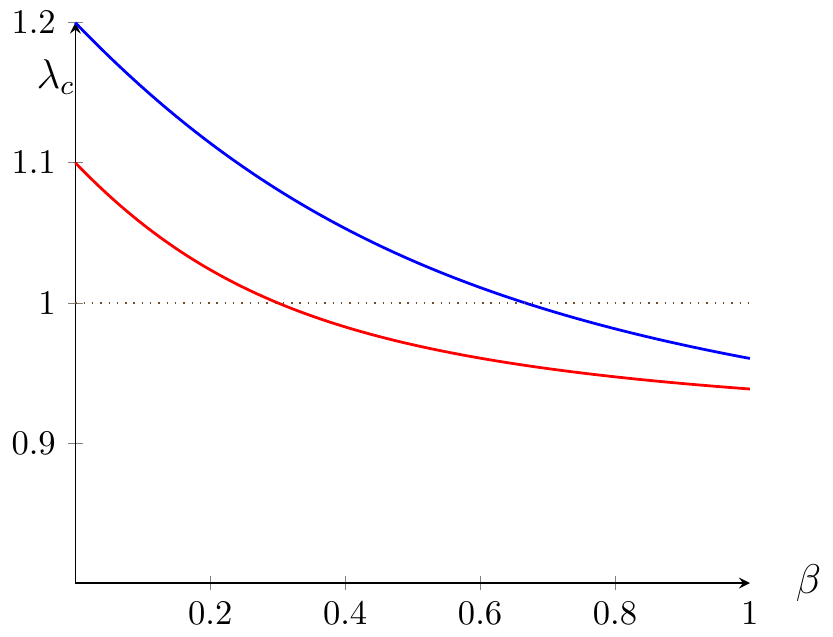}}}
\centering \mbox{\subfigure[Simulated size of the largest component.]{\includegraphics[height=5cm]{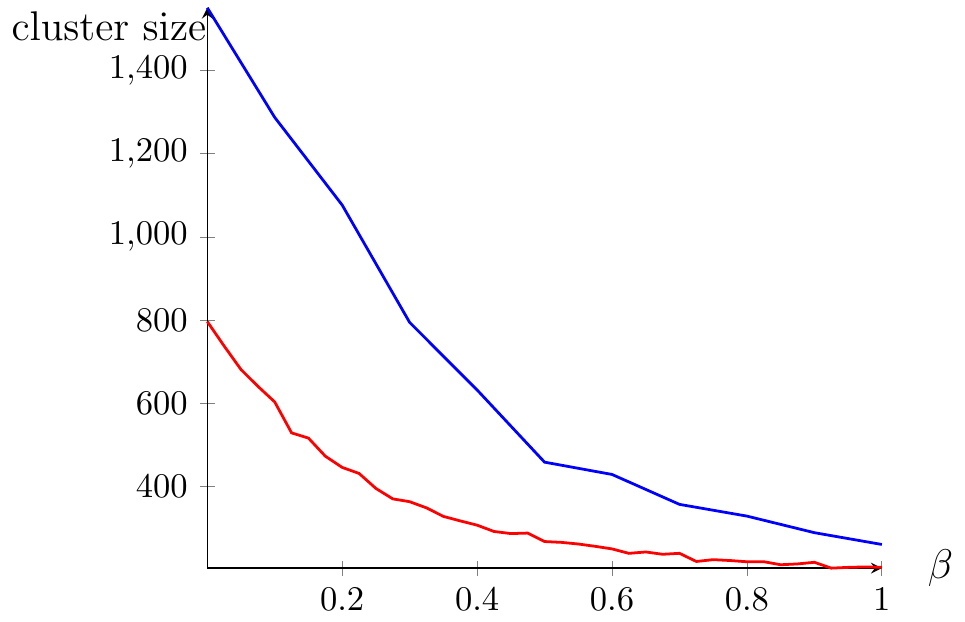}}}
\caption{The critical parameter $\lambda_c$ and the simulated size of the largest component plotted against $\beta$ for the two-type Erd\H{o}s-R\'{e}nyi graph with $p_1=p_2=1/2$ and $\mu_1=0.5$, $\mu_2=1.2$ (blue) and $\mu_1=0.7$, $\mu_2=1.1$ (red).}\label{fig:ER_Ex2}
\end{figure}

\begin{figure}[p!]
\centering \mbox{\subfigure[The critical parameter.]{\includegraphics[height=5cm]{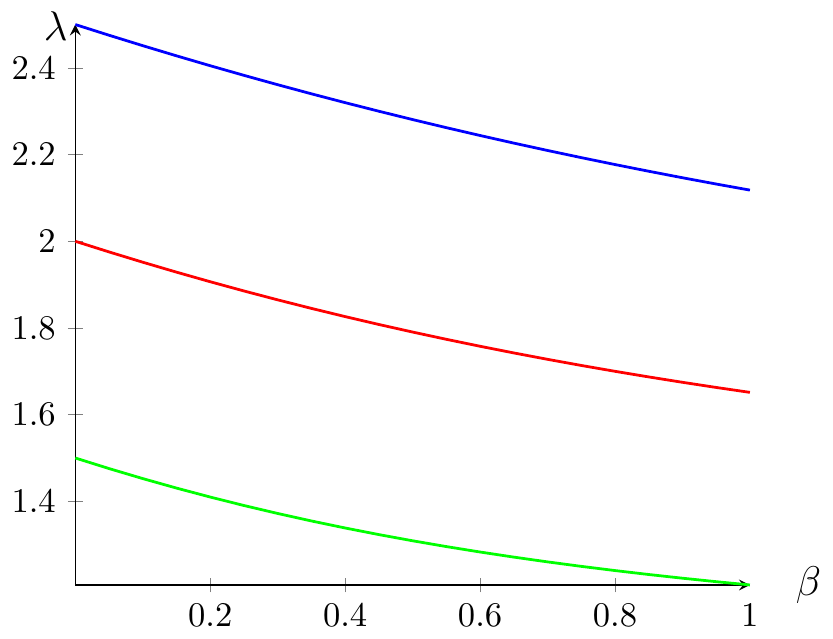}}}
\centering \mbox{\subfigure[Simulated size of the largest component.]{\includegraphics[height=5cm]{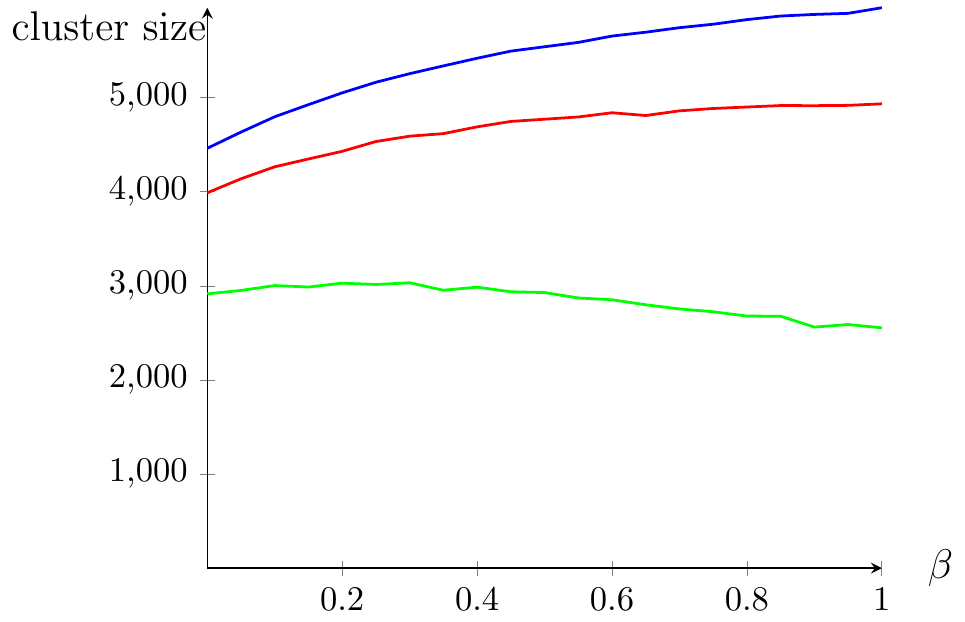}}}
\caption{The critical parameter $\lambda_c$ and the simulated size of the largest component plotted against $\beta$ for the two-type Erd\H{o}s-R\'{e}nyi graph with $p_1=p_2=1/2$ and $\mu_1=0.5$ for $\mu_2=2.5$ (blue) $\mu_2=2$ (red) and $\mu_2=1.5$ (green).}\label{fig:ER_Ex2(2)}
\end{figure}

\noindent \textbf{Example 2.2: One subcritical and one supercritical type.} Again let $p_1=1/2$, but take $\alpha_1$ and $\alpha_2$ to be (potentially) different. The asymptotic average degree for type $i$ vertices is $\mu_i=(\alpha_i+\beta)/2$. First consider the case when type 1 is subcritical ($\mu_1<1$) and type 2 is supercritical ($\mu_2>1$). Figure \ref{fig:ER_Ex2} shows plots of $\lambda_c$ and simulated sizes of the largest component against $\beta$ for two different combinations of expected degrees: $\mu_1=0.5$, $\mu_2=1.2$ and $\mu_1=0.7,\mu_2=1.1$. We see that $\lambda_c$ decreases as $\beta$ increases and eventually drops below 1, and the component size decreases as $\beta$ increases. This is because the subcritical type 1 part of the population restrains the type 2 part when more edges are used to connect vertices of different type. Figure \ref{fig:ER_Ex2(2)} contains analogous pictures for the case when $\mu_1=0.5$ and $\mu_2\in\{1.5,2,2.5\}$. The critical parameter $\lambda_c$ is still decreasing in $\beta$, for the same reason, but here the influence of type 1 is not sufficient to bring $\lambda_c$ below 1. The size of the largest component is now slightly increasing in $\beta$ for the two larger values of $\mu_2$. This is because superfluous edges connecting type 2 vertices can be shifted to include type 1 vertices in the giant component without type 2 vertices being disconnected. When $\mu_2=1.5$, the surplus of type 2 edges is not large enough to allow for connections to the subcritical type 1, and the component size decreases with $\beta$.\hfill$\Box$

\begin{figure}[p!]
\centering \mbox{\subfigure[The critical parameter.]{\includegraphics[height=5cm]{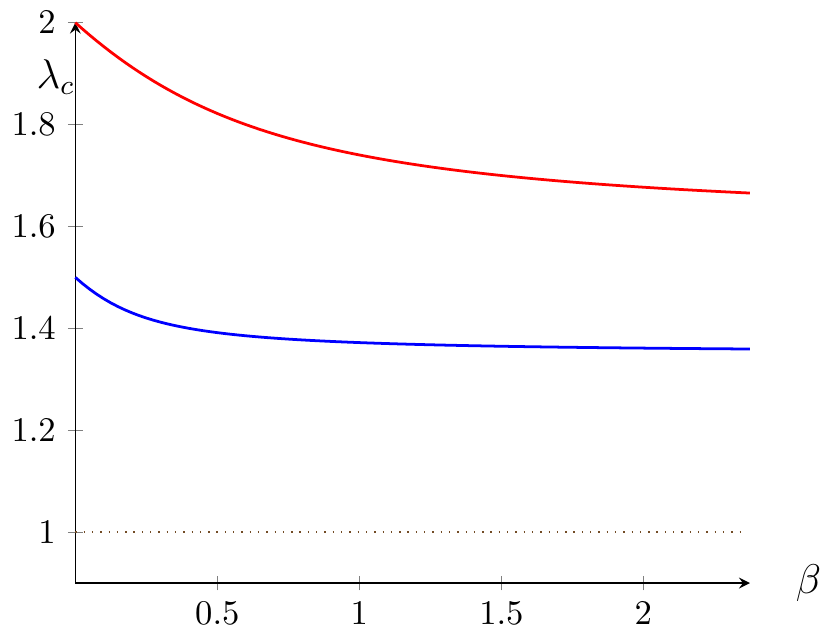}}}
\centering \mbox{\subfigure[Simulated size of the largest component.]{\includegraphics[height=5cm]{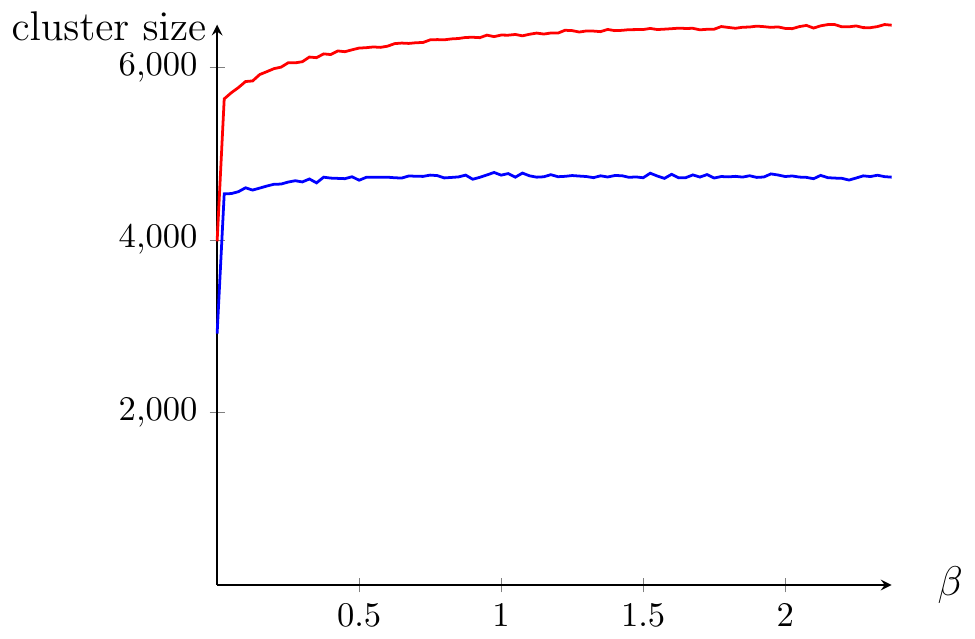}}}
\caption{The critical parameter $\lambda_c$ and the simulated size of the largest component plotted against $\beta$ for the two-type Erd\H{o}s-R\'{e}nyi graph with $p_1=p_2=1/2$ and $\mu_1=1.2$ for $\mu_2=1.5$ (blue) and $\mu_2=2.0$ (red).}\label{fig:ER_Ex3}
\end{figure}

\noindent \textbf{Example 2.3: Two supercritical types.} Next we consider the case when both types are supercritical. Figure \ref{fig:ER_Ex3} shows pictures for two cases: $\mu_1=1.2$, $\mu_2=1.5$ and $\mu_1=1.2$, $\mu_2=2$ (with $p_1=p_2=1/2$). The critical parameter is again decreasing in $\beta$ because of the larger influence of the type with smaller expected degree. Here, the decrease occurs mainly for smaller values of $\beta$, when the types start to mix, and the curve then flattens out. As for the component size, when both types are supercritical, the giant type 1 and type 2 components merge as soon as $\beta$ becomes positive, which leads to an increase in the size of the largest component. The size of the largest component then does not change much as $\beta$ increases.\hfill$\Box$

\begin{figure}[p!]
\centering \mbox{\subfigure[The critical parameter.]{\includegraphics[height=5cm]{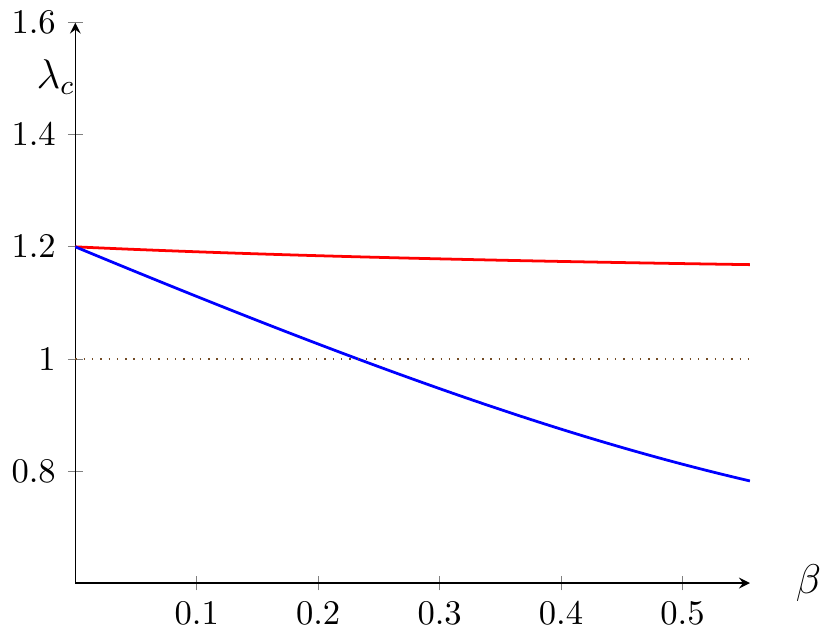}}}
\centering \mbox{\subfigure[Simulated size of the largest component.]{\includegraphics[height=5cm]{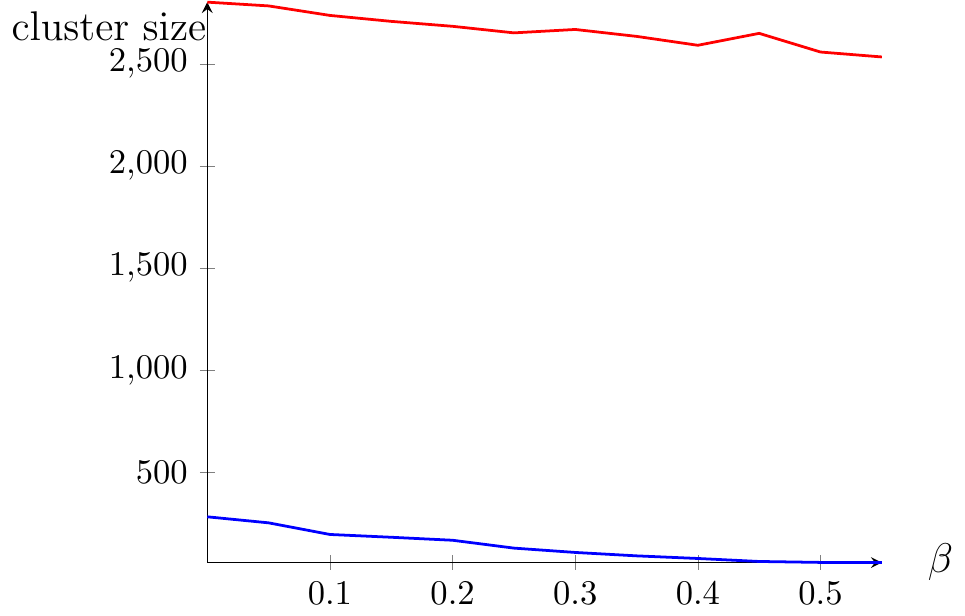}}}
\caption{The critical parameter $\lambda_c$ and the simulated size of the largest component plotted against $\beta$ for the two-type Erd\H{o}s-R\'{e}nyi graph with $\mu_1=0.5$ and $\mu_2=1.2$ for $p_1=0.1$ (blue) and $p_1=0.9$ (red).}\label{fig:ER_Ex4}
\end{figure}

\noindent \textbf{Example 2.4: Minorities.} Next we consider a situation where one group is a minority. We fix $\mu_1=0.5$ and $\mu_2=1.2$. Figure \ref{fig:ER_Ex4}(a) shows plots of $\lambda_c$, given in \eqref{eq:ERl_allm}, against $\beta$ for $p_1=0.9$ and $p_1=0.1$, respectively. We see that, when the subcritical type 1 part occupies the majority of the graph, then $\lambda_c$ drops below 1 faster compared to the symmetric case (cf.\ Figure \ref{fig:ER_Ex2}(a)) while, when the proportions are interchanged so that type 1 is the minority, then the influence of type 1 is not sufficient to make the graph subcritical. Simulated sizes of the largest component are plotted in Figure \ref{fig:ER_Ex4}(b). Note that, for $p_1=0.9$, the expected number of type 2 vertices is only 1000 and, for $\beta=0$, about half of these create a large component, whose size then drops to 0 as $\beta$ increases. For $p_1=0.1$, type 2 dominates the graph and the relative size of the largest component does not change much as $\beta$ increases.\hfill$\Box$

\begin{figure}
\centering \mbox{\subfigure[Critical parameter.]{\includegraphics[height=5cm]{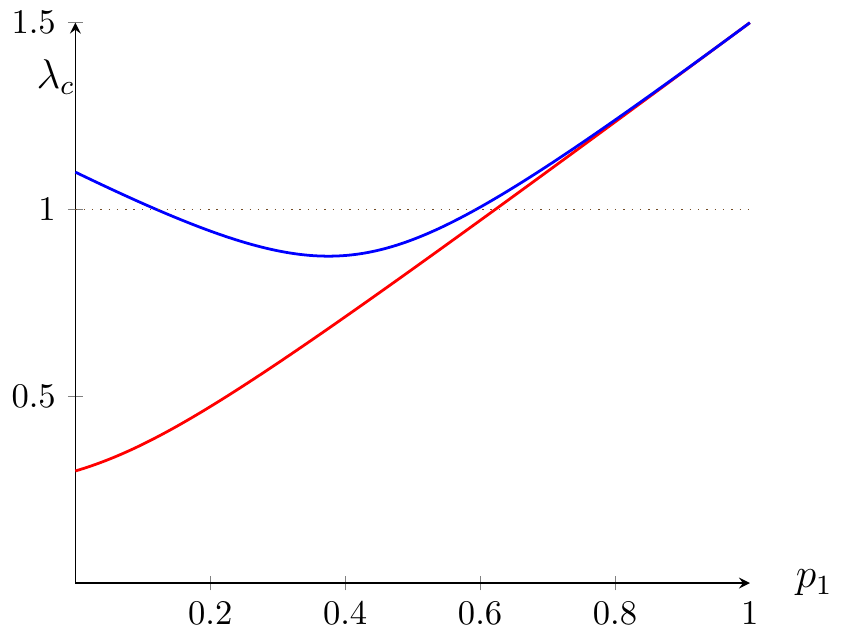}}}\par
\centering \mbox{\subfigure[Simulated size of largest component.]{\includegraphics[height=5cm]{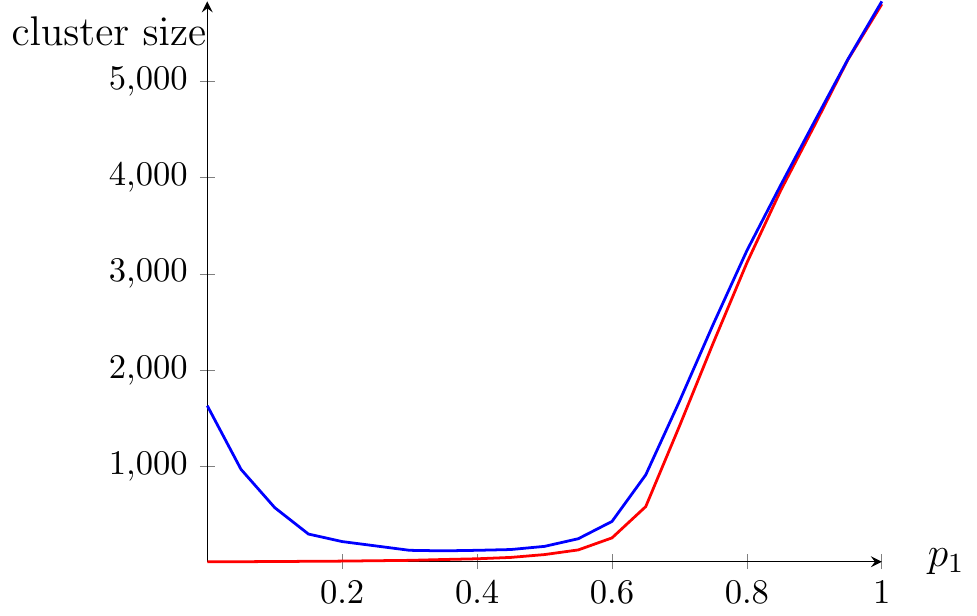}}}
\centering \mbox{\subfigure[Simulated size of 2nd (solid lines) and 3rd (dashed lines) largest component.]{\includegraphics[height=5cm]{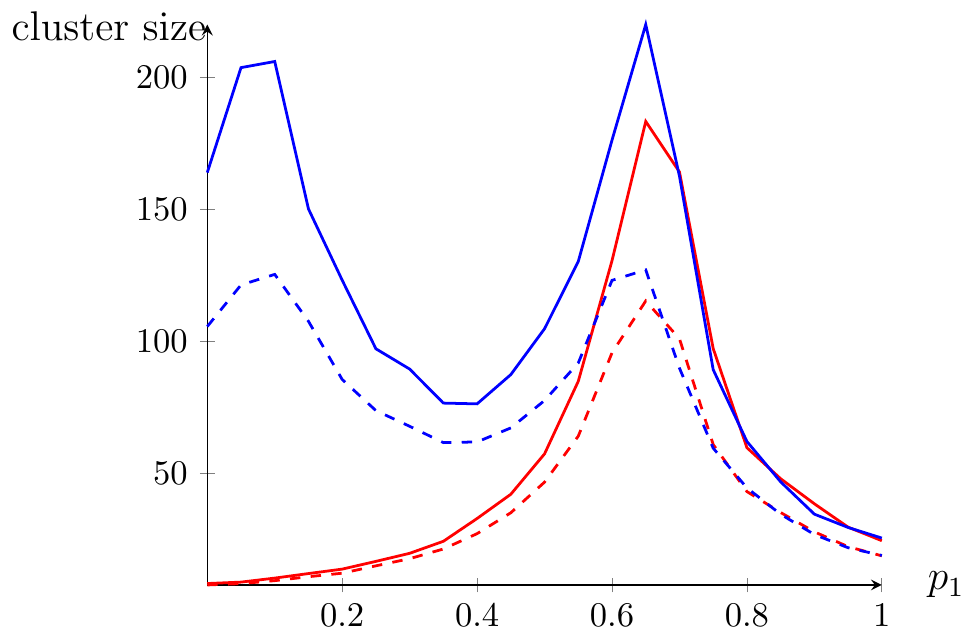}}}
\caption{The critical parameter $\lambda_c$ and simulated component sizes plotted against $p_1$ for the two-type Erd\H{o}s-R\'{e}nyi graph with $\beta=0.5,\alpha_2=1.5$ and $\alpha_1=1.1$ (blue), $\alpha_1=0.3$ (red).}\label{fig:ER_Ex5}
\end{figure}

\noindent \textbf{Example 2.5: Varying $p_1$.} Finally, we consider a different type of example where we, instead of varying $\beta$, vary the proportions of the types, while keeping $\alpha_1$, $\alpha_2$ and $\beta$ fixed. This means that, as the proportions change, so do the average degrees of the types. In Figure \ref{fig:ER_Ex5}(a), the critical parameter $\lambda_c$ is plotted against $p_1$ for $\alpha_1=1.5,\alpha_2=1.1,\beta=0.5$ (blue) and $\alpha_1=1.5,\alpha_2=0.3,\beta=0.5$ (red).

In the first case, both types prefer to connect to their own type, and type 1 is more social than type 2. As $p_1$ grows, the average type 1 degree increases from $\beta$ to $\alpha_1$ while the average type 2 degree decreases from $\alpha_2$ to $\beta$. We see that the critical parameter first decreases (until the average degrees coincide) and then instead increases. In particular, there is an interval for $p_1$ where the graph is subcritical -- the previously dominating type 2 suffers from a shortage of other type 2 vertices to connect to, and this has not yet been compensated for by the influx of the more social type 1.

In the second case, both types have a larger probability of connecting to type 1, so both average degrees increase with $p_1$ and hence so does $\lambda_c$. Simulated sizes of the largest component are plotted against $p_1$ in Figure \ref{fig:ER_Ex5}(b) and behave as expected in relation to the critical parameter. We note that, for large values of $p_1$, the difference between the two cases is small, since the graph is then in both cases largely controlled by type 1.

In Figure \ref{fig:ER_Ex5}(c), we have also included a plot of the size of the second and third largest component against $p_1$ (although it is not an ambition of this work to give a systematic analysis of the smaller order components). For the standard Erd\H{o}s-R\'{e}nyi graph it is well-know that the sizes of the smaller order components are of order $n^{2/3}$ when the graph is critical and of order $\log n$ otherwise, that is, the smaller order components peak at criticality. The plot nicely illustrates that this is the case also for the two-type version.\hfill$\Box$

\section{The configuration model}

We continue with the two-type configuration model, and first recall the definition of the model: A vertex is type $i$ ($i=1,2$) with probability $p_i$, and type $i$ vertices are independently assigned half-edges according to a given degree distribution $F_i$ with finite mean. Half-edges at type $i$ vertices are independently labeled $i$ or $i'$ with probability $\xi_i$ and $1-\xi_i$, respectively. Half-edges with label $i$ are then paired to each other uniformly at random, and half-edges with label $1'$ are paired to half-edges with label $2'$ uniformly at random.

The pairing procedure might leave a number of half-edges unpaired. Specifically, at most one half-edge with label 1 (2) can be left without a partner (this happens when the number of label 1 (2) half-edges is odd), and pairing half-edges with label 1' and 2' to each other will leave $|L_{1'}-L_{2'}|$ half-edges unpaired, where $L_{i'}$ denotes the total number of half-edges with label $i'$. We erase all un-paired half-edges. This might of course affect the degree distributions. The change coming from erasing a single stub with label 1 (2) however is clearly negligible in the limit as $n\to\infty$. Let $D_i^{\sss \rm{tot}}$ denote the total degree of the type $i$ vertices (before erasing) and write $\mu_i$ for the mean of $F_i$. We have $L_{i'}\sim \mbox{Bin}\left(D_i^{\sss\rm{tot}}, 1-\xi_i\right)$ so that $\E[L_{i'}]=np_i\mu_i(1-\xi_i)$. In order to control the effects of erasing half-edges with label 1' or 2', we require that $\E[L_{1'}]=\E[L_{2'}]$, that is,
\begin{equation}\label{conf_balance}
p_1\mu_1(1-\xi_1)=p_2\mu_2(1-\xi_2).
\end{equation}
The number of vertices that have at least one un-paired half-edge -- and whose degrees are thereby affected by the erasing procedure -- is dominated from above by $|L_{1'}-L_{2'}|$. Assuming \eqref{conf_balance} it is not hard to see that the fraction of vertices that have at least one un-paired half-edge converges to 0 in probability and in mean (essentially this follows from the law of large numbers). We conclude that \eqref{conf_balance} ensures that the distribution of the degree of a randomly chosen vertex of a given type is asymptotically the same as the input distribution. Note that the model can produce self-loops and multiple edges. We allow for this, but remark that, when the degrees have finite mean, self-loops and multiple edges can also be erased without affecting the degree distributions; see \cite[Section 7.3]{RemcoI}.

Just like the two-type Erd\H{o}s-R\'{e}nyi graph, the model undergoes a phase transition and again the critical parameter can be obtained from theory for multi-type branching processes. We identify the components in the graph by picking a vertex at random and then exploring its component in generations -- first the neighbors are explored, then the neighbors of the neighbors, and so on. The procedure is similar to the derivation of the phase-transision in the one-type configuration model, where the early stages of the exploration can be approximated by a one-type branching process; see e.g.\ \cite[Section 4]{RemcoII} for details. In our case, as $n\to\infty$, the early stages of the exploration are well approximated by a two-type branching process with the following offspring distributions. The initial individual has offspring distribution prescribed by the degree and type distributions in the model. In the second and later generations, the offspring distributions become size biased: if the parent is type $i$, then the total number of children -- corresponding to neighboring vertices that have not yet been explored -- is distributed as $\widetilde{D}_i-1$, where $\widetilde{D}_i$ is a size biased version of $D_i\sim F_i$, that is,
$$
\PP(\widetilde{D}_i=d)=\frac{d\PP(D_i=d)}{\mu_i}.
$$
To see this note that, by construction of the graph, the neighbors of a given vertex are determined in that half-edges are chosen uniformly at random from the set of all half-edges. Hence the degree of an individual in the second and later generation is given by the degree of the vertex of a randomly chosen half-edge. The probability of connecting to a vertex with degree $d$ is then proportional to $d$, explaining the size biasing effect; see \cite[Section 4]{RemcoII} for more details.

We have that $\E[\widetilde{D}_i-1]=(\E[D_i^2]-\mu_i)/\mu_i=:\nu_i$, which is assumed to be finite for both types in what follows; see below for comments and simulations of cases where $\nu_i=\infty$. Furthermore, a given child is type $i$ with probability $\xi_i$ and type $i^c$ with probability $1-\xi_i$ (recall the notation $i^c$ from the previous section). Let $m_{ik}$ denote the expected number of type $k$ children of a type $i$ individual, and write $M=\{m_{ik}\}$ for the reproduction matrix. Then,
$$
M=\begin{pmatrix}
\xi_1\nu_1 & (1-\xi_1)\nu_1  \\
(1-\xi_2)\nu_2 & \xi_2\nu_2
\end{pmatrix}.
$$
The branching process has a strictly positive probability of exploding if and only if the largest positive eigenvalue of $M$ is strictly larger than 1. This eigenvalue corresponds to the critical parameter for the graph model: if the eigenvalue is strictly larger than 1, then the graph will contain a unique giant component with size of order $n$ (hitting this giant component in the exploration of the graph corresponds to explosion in the approximating branching process), while if the eigenvalue is smaller than or equal to 1, then all components are sublinear in $n$. We anticipate that this link can be made fully rigorous by minor modifications of the methods in \cite{Ball_Sirl}. Note that the assumption \eqref{conf_balance} is necessary for the matrix $M$ to provide a good approximation of the exploration of the graph. We conclude, after standard calculations, that the critical parameter for the model is given by
$$
\lambda_c=\frac{\xi_1\nu_1+\xi_2\nu_2}{2}+\sqrt{\left(\frac{\xi_1\nu_1+\xi_2\nu_2}{2}\right)^2+\nu_1\nu_2(1-\xi_1-\xi_2)}.
$$
If one or both degree distributions has infinite second moment, so that $\nu_i=\infty$ for $i=1$ or $i=2$, then the graph contains a giant component: The one-type configuration model generates a giant component for degree distributions with infinite second moment; see \cite[Section 4]{RemcoII}. Indeed, the exploration process can then be approximated by a branching process with infinite mean, which has a positive probability of exploding. If, say, $\nu_1=\infty$ in the two-type model, then clearly the number of half-edges with label 1 at a type 1 vertex has infinite second moment. Since the restriction of the two-type graph to type 1 vertices and edges between them has the same distribution as a one-type graph with degree distribution given by the number of half-edges with label 1 at a type 1 vertex, it follows from the result for the one-type model that the two-type model will generate a giant component containing a positive fraction of the type 1 vertices (and thereby also a positive fraction of the type 2 vertices as soon as $\xi_2\neq 1$).

The limiting size of the giant component can again be expressed in terms of equations involving generating functions for the approximating branching process, but also in this case we prefer simulations.

The extent to which type $i$ vertices connect to other type $i$ vertices is controlled by the parameter $\xi_i$ -- a larger value of $\xi_i$ means that edges incident to a type $i$ vertex are more likely to connect to other type $i$ vertices. As for the Erd\H{o}s-R\'{e}nyi graph, we keep the degree distributions and proportions of the types fixed in the examples. This means that the expected number of edges in the graph is fixed. We then vary the extent to which edges connect vertices of similar/different type by varying e.g.\ $\xi_1$ and adjust $\xi_2$ accordingly to ensure \eqref{conf_balance}. Thus, when $\xi_1$ increases, also $\xi_2$ increases. We then investigate how the critical parameter and the size of the giant component are affected.

\noindent \textbf{Example 3.1 ($\nu_1=\nu_2$).} We begin by observing that, for distributions $F_1$ and $F_2$ with $\nu_1=\nu_2=\nu$, using \eqref{conf_balance} and some algebra, one can deduce that $\lambda_c=\nu$. In this case, the threshold is hence not affected by the value of $\xi_1$, but coincides with the threshold for the standard configuration model.\hfill$\Box$

\begin{figure}
\centering \mbox{\subfigure[The critical parameter.]{\includegraphics[height=5cm]{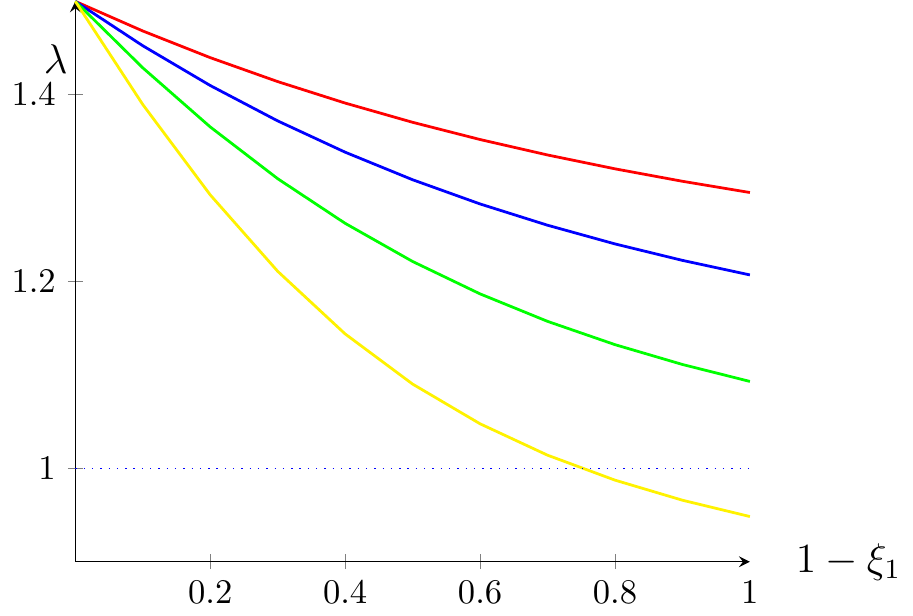}}}
\centering \mbox{\subfigure[Simulated size of the largest component.]{\includegraphics[height=5cm]{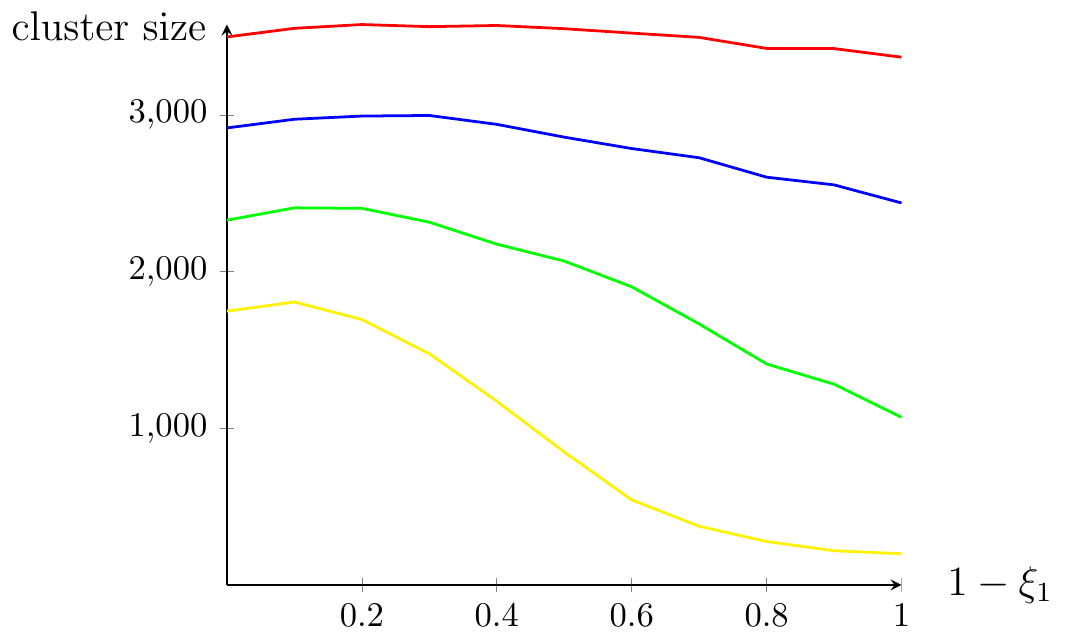}}}
\caption{The critical parameter $\lambda_c$ for the two-type configuration model with $\mu_1=\nu_1=0.5$ and $\mu_2=\nu_2=1.5$ plotted against $1-\xi_1$ for varying proportions of the types: $p_1=0.4$ (red), $p_1=0.5$ (blue), $p_1=0.6$ (green), $p_1=0.7$ (yellow). Simulated component sizes for $F_1=\rm{Poisson}(0.5)$ and $F_2=\rm{Poisson}(1.5)$.}\label{fig:CM_popo_subsup}
\end{figure}

\noindent \textbf{Example 3.2 ($\nu_1\neq \nu_2$).} Next we consider the case when $\nu_1\neq \nu_2$. In Figure \ref{fig:CM_popo_subsup}(a), $\lambda_c$ is plotted against $1-\xi_1$ for distributions with $\nu_1=\mu_1=0.5$ and $\nu_2=\mu_2=1.5$ for different values of the type 1 proportion (we plot against $1-\xi_1$ rather than $\xi_1$ to make the plots more analogous to the ones for the Erd\H{o}s-R\'{e}nyi graph). We see that $\lambda_c$ decreases when the types become more inclined to attach to each other. Just as in the Erd\"{o}s-Renyi graph, the reason is that the influence of a part of the population with a smaller average degree gradually increases. For $p_1=0.7$, the threshold $\lambda_c$ drops below 1 when $1-\xi_1$ is around 0.75, while for the smaller values of $p_1$, the subcritical type 1 is not able to bring $\lambda_c$ below 1.\hfill$\Box$

The critical parameter depends on the degree distributions only through their means and second moments. The size of the largest component however may be affected also by other properties of the degree distribution. We have simulated component sizes for a few degree distributions with different tail behavior. Here we have also included distributions that do not have a finite second moment, so that $\nu=\infty$. Throughout, graphs of size $n=10,000$ are simulated and the pictures are based on averages of 100 iterations per parameter value; see \cite{code} for a link to the code.

We study combinations of Poisson distributions and power-law distributions of Yule-Simon type. Power-laws are often observed in applications and it is therefore of interest to study this distribution type. A random variable $D$ has a Yule-Simon distribution if $\PP(D=k)=\rho B(k,\rho+1)$ for $k=1,2,\ldots$ and some $\rho>0$, where $B(\cdot,\cdot)$ is the beta function. As $k\to\infty$, the probability $\PP(D=k)$ decays like $k^{-(\rho+1)}$. We have $\E[D]=\rho/(\rho-1)$ for $\rho>1$, so that a larger tail exponent hence implies a smaller mean. Values $\E[D]\in(2,\infty)$ correspond to $\rho\in(1,2)$, which implies infinite variance, while $\E[D]\in(1,2)$ correspond to $\rho\in(2,\infty)$, implying finite variance. For $\rho>2$, the critical parameter is given by $\nu=2/(\rho-2)$, and the graph is subcritical if $\nu<1$, which translates into $\rho>4$, that is, $\E[D]<4/3$.

\begin{figure}
\centering \mbox{\subfigure[$\mu_1=2.5$, $\mu_2=1.2$ (blue)\newline\hspace*{0.5cm} $\mu_1=2.0$, $\mu_2=1.2$ (red)\newline\hspace*{0.5cm} $\mu_1=1.5$, $\mu_2=1.2$ (green)]{\includegraphics[height=4.7cm]{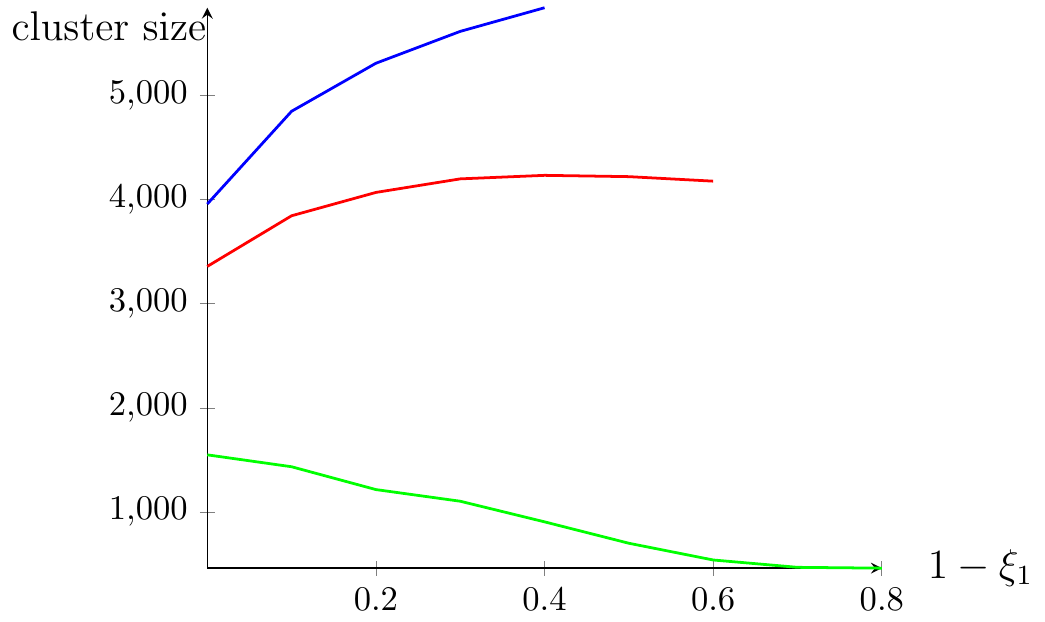}}}
\centering \mbox{\subfigure[$\mu_1=2.0$, $\mu_2=2.5$ (blue)\newline\hspace*{0.5cm} $\mu_1=1.5$, $\mu_2=2.0$ (red)\newline\hspace*{0.5cm} $\mu_1=1.5$, $\mu_2=2.5$ (green)]{\includegraphics[height=4.7cm]{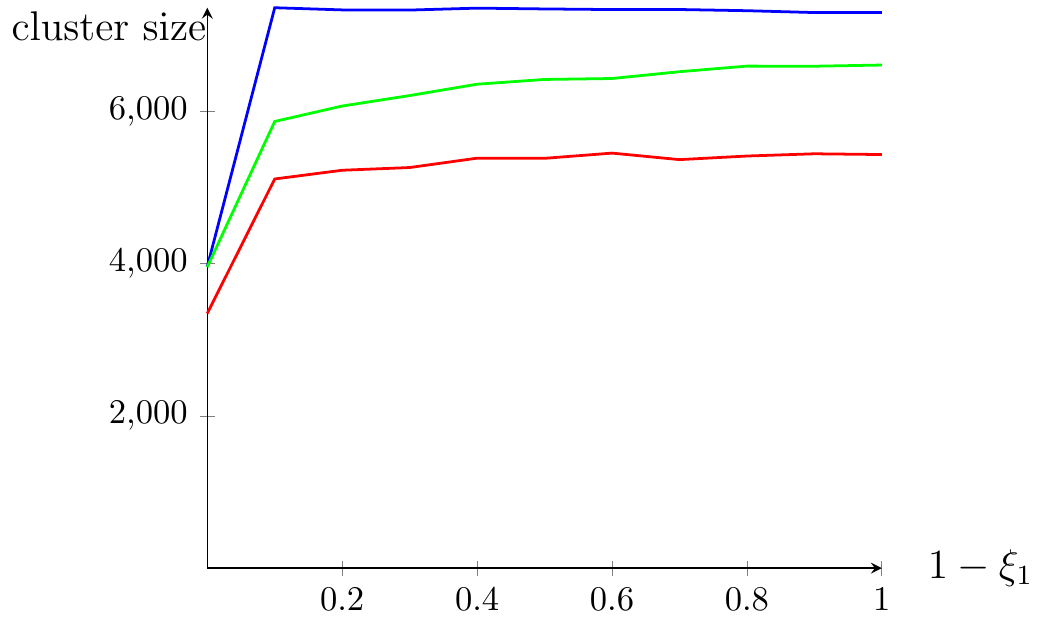}}}
\caption{Simulated component sizes plotted against $1-\xi_1$ in the two-type configuration model with $p_1=p_2=1/2$ and $F_1=\mbox{Yule-Simon}(\mu_1)$ and $F_2=\mbox{Yule-Simon}(\mu_2)$. In (a), $F_1$ is subcritical and $F_2$ supercritical, in (b) both are supercritical.}\label{fig:CM_YSYS}
\end{figure}

\begin{figure}
\centering \mbox{\subfigure[$\mu_1=2.5$, $\mu_2=1.1$ (blue)\newline\hspace*{0.5cm} $\mu_1=2.0$, $\mu_2=1.1$ (red)\newline\hspace*{0.5cm} $\mu_1=1.5$, $\mu_2=1.1$ (green)\newline\hspace*{0.5cm} $\mu_1=1.2$, $\mu_2=1.1$ (yellow)]{\includegraphics[height=4.7cm]{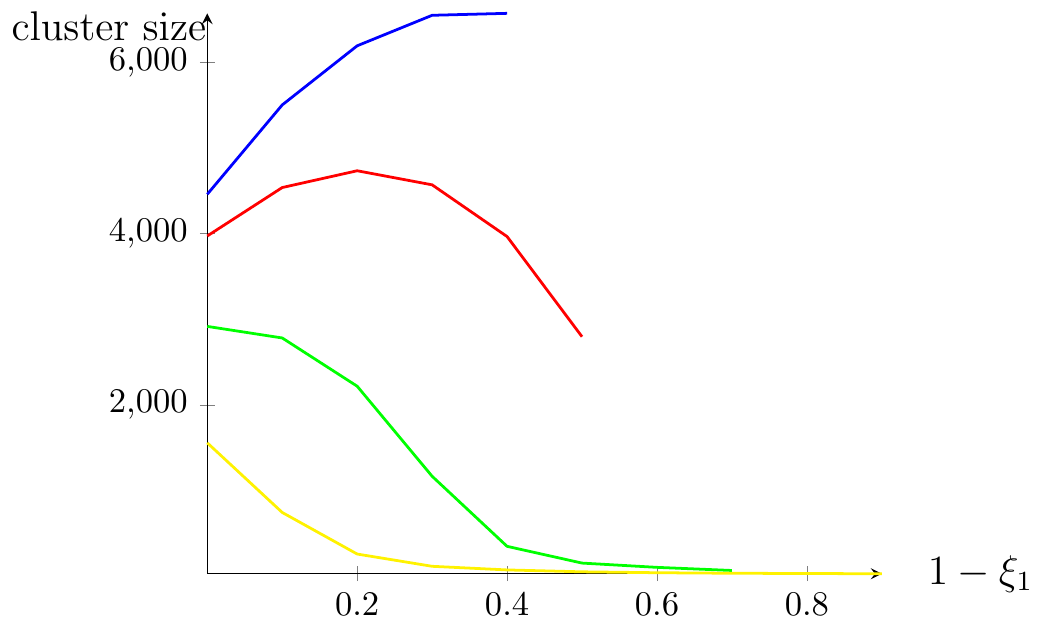}}}
\centering \mbox{\subfigure[$\mu_1=2.0$, $\mu_2=2.5$ (blue)\newline\hspace*{0.5cm} $\mu_1=1.5$, $\mu_2=2.5$  (red)\newline\hspace*{0.5cm} $\mu_1=1.2$, $\mu_2=2.5$ (green)]{\includegraphics[height=4.7cm]{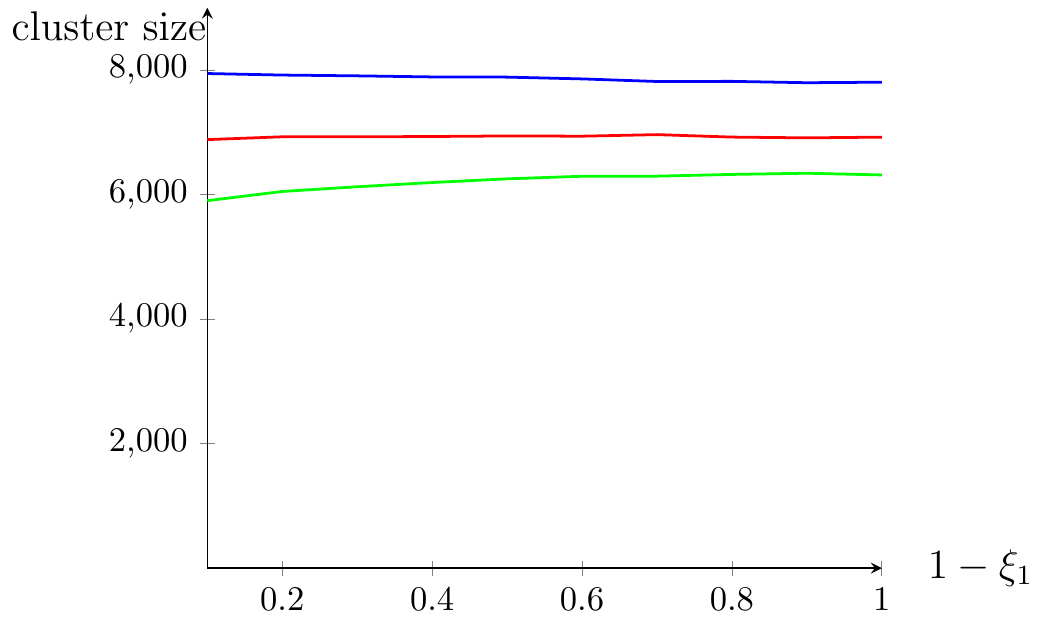}}}
\caption{Simulated component sizes plotted against $1-\xi_1$ in the two-type configuration model with $p_1=p_2=1/2$ and $F_1=\mbox{Poisson}(\mu_1)$ and $F_2=\mbox{Yule-Simon}(\mu_2)$. In (a), the Poisson is supercritical and the Yule-Simon is subcritical, in (b) both are supercritical.}\label{fig:CM_PoYS}
\end{figure}

Figure \ref{fig:CM_popo_subsup}(b) shows final sizes plotted against $1-\xi_1$ for a few different values of $p_1$ when $F_1=\rm{Poisson}(0.5)$ and $F_2=\rm{Poisson}(1.5)$, that is, one subcritical and one supercritical Poisson distribution. We see that the component size decreases as $1-\xi_1$ increases (so that more edges connect vertices of different types) and that the decrease is larger for larger values of $p_1$ (when the subcritical type 1 is more dominant). The reason is the same as for the Erd\H{o}s-R\'{e}nyi graph: attaching to the subcritical type means that vertices of the supercritical type are cut off from the largest component.

Figure \ref{fig:CM_YSYS} shows final sizes plotted against $1-\xi_1$ for $p_1=1/2$ when both $F_1$ and $F_2$ are Yule-Simon distributions. In (a), $F_1$ is supercritical with varying mean $\mu_2$ and $F_2$ is subcritical with mean $\mu_1=1.2$ (in some cases the plots are not continued for large values of $1-\xi_1$ since it is then not possible to define $\xi_2$ as to satisfy \eqref{conf_balance}). We see that the component size is increasing for large values of $\mu_1$ and decreasing for a smaller value of $\mu_1$. In (b), both types are supercritical and, after the initial increase when the two giant components merge (delayed in the picture because of our simulation grid), the component size is essentially constant. Both these findings are analogous to the Erd\H{o}s-R\'{e}nyi case.

Figure \ref{fig:CM_PoYS}, shows component sizes for $p_1=1/2$ when $F_1$ is Poisson and $F_2$ a Yule-Simon distribution. In (a), the Poisson is supercritical with varying mean $\mu_1$ and the Yule-Simon is subcritical with mean $\mu_2=1.1$. We see that, as $\mu_1$ decreases, the largest component naturally becomes smaller. Interestingly, when $\mu_1=2.0$, the component size is not monotone in $1-\xi_1$, but goes from being increasing to decreasing: when the mixing between the types becomes too large, the subcritical type 1 starts to cut off parts of the giant component. In (b), both distributions are supercritical, leaving the component size roughly constant as $1-\xi_1$ varies.

\section{Preferential attachment}

Finally, we consider the two-type preferential attachment model. Starting at time $t=1$ with one vertex of each type connected by an edge, a new vertex $v_t$ with one edge on it arrives in the graph at each integer time $t\geq 2$. The new vertex is type $i$ with probability $p_i$ ($i=1,2$) and then connects to another type $i$ vertex with probability $\theta_i$ and to a type $i^c$ vertex with probability $1-\theta_i$. When a new vertex has decided which type to connect to, it chooses a vertex of that type proportionally to degree. Let $T_t$ denote the type that vertex $v_t$ chooses to connect to, write $G(t)$ for the graph at time $t$ and $D_s(t)$ for the degree of vertex $v_s$ in $G(t)$. Furthermore, let $L_i(t)$ denote the total degree of all type $i$ vertices in $G(t)$. Then,
$$
\PP(v_{t+1}\to v_s|G(t), T_{t+1}=i)= \left\{ \begin{array}{ll}
                      \frac{D_s(t)}{L_i(t)} & \mbox{if $v_s$ is type $i$};\\
                      0 & \mbox{otherwise}.
                    \end{array}
            \right.
$$
Note that we obtain a connected graph for all parameter values.

The large interest in preferential attachment models during the last decade stems from the fact that they typically give rise to power-law degree sequences and thereby reproduce an important feature of many empirical networks. Below, we heuristically derive expressions for the expected degree sequences of the types in our two-type version, confirm that they are indeed power-laws and identify the exponents. We have not attempted to do a rigorous analysis and prove concentration. This may be possible using traditional martingale methods, or an approach based on general branching processes \cite{D,RT}.

Write $N^{\sss(k)}_i(t)$ for the number of type $i$ vertices with degree $k$ at time $t$. How does $N^{\sss(k)}_i(t)$ change when $v_{t+1}$ arrives? It increases by 1 if $v_{t+1}$ attaches to a type $i$ vertex with degree $k-1$ (and if $v_{t+1}$ is type $i$ and $k=1$) and decreases by 1 if $v_{t+1}$ attaches to a type $i$ vertex with degree $k$. The probability that $v_{t+1}$ attaches to a type $i$ vertex is $a_i:=p_i\theta_i+p_{i^c}(1-\theta_{i^c})$. In what follows we shall assume that $a_i>0$. We obtain that
\begin{equation}\label{mean_change}
\E[N^{\sss(k)}_i(t+1)|G(t)]=N^{\sss(k)}_i(t)+\frac{a_i(k-1)N^{\sss(k-1)}_i(t)}{L_i(t)}
-\frac{a_ikN^{\sss(k)}_i(t)}{L_i(t)}+p_i\mathbf{1}_{\{k=1\}}.
\end{equation}
The total degree of type $i$ vertices increases by 1 if a new vertex is type $i$ and it also increases by 1 if the new vertex attaches to an existing type $i$ vertex. The expected change in each step is hence $p_i+a_i$ and it follows from the law of large numbers that $L_i(t)/t\to p_i+a_i$. Write $N_i(t)$ for the total number of type $i$ vertices at time $t$. Clearly $N_i(t)/t\to p_i$.

Now assume that $N^{\sss(k)}_i(t)/t$ converges in mean to some limit $p_ir_i^{\sss (k)}$, where $r_i^{\sss (k)}$ represents the asymptotic expected fraction of the type $i$ vertices that has degree $k$. In fact, we will assume that the convergence holds also in probability\footnote{Both these assumptions are natural in view of previous work on preferential attachment models. If $N_i(t)$ and $L_i(t)$ are concentrated around their means, we have that $\E[N^{\sss(k)}_i(t)]/t$ satisfies a recursion obtained from (\ref{mean_change}). It is then reasonable to believe that the solution of the recursion, which we denote by $p_ir_i^{\sss (k)}$, gives a good approximation of $\E[N^{\sss(k)}_i(t)]/t$ for large $t$ so that, in particular, $\E[N^{\sss(k)}_i(t)]/t\to p_ir_i^{\sss (k)}$. Concentration results ensuring convergence in probability are traditionally established using margingale methods first appearing in \cite{BRST}.}. We then obtain that
$$
\frac{N^{\sss(k)}_i(t)}{L_i(t)}\stackrel{p}{\to} \frac{p_ir_i^{\sss(k)}}{p_i+a_i}.
$$
Write $b_i=p_i+a_i$. We now take expectation on both sides in \eqref{mean_change} and let $t\to\infty$. Using that $\E[N^{\sss(k)}_i(t)]$ is close to $tp_ir_i^{\sss(k)}$ when $t$ is large, and also using bounded convergence, we obtain that
$$
r_i^{\sss(k)}=\frac{(k-1)a_i}{b_i}r_i^{\sss(k-1)}-\frac{ka_i}{b_i}r_i^{\sss(k)}+
\mathbf{1}_{\{k=1\}},
$$
which can be rewritten as
$$
r_i^{\sss(k)}=\frac{k-1}{b_i/a_i+k}r_i^{\sss(k-1)}+
\frac{1}{1+ka_i/b_i}\mathbf{1}_{\{k=1\}}.
$$
Iterating this recursion and recalling the definition and properties of the Gamma function (see e.g.\ \cite[Chapter 6]{functions}) yields that
$$
r_i^{\sss(k)}=\frac{1}{1+a_i/b_i}\prod_{j=1}^{k-1}\frac{j}{b_i/a_i+1+j}=\frac{1}{1+a_i/b_i}\cdot
\frac{\Gamma(k)\Gamma(b_i/a_i+2)}{\Gamma(1)\Gamma(k+1+b_i/a_i)}.
$$
We have hence arrived at an explicit expression for the asymptotic expected degree distributions for the types in the graph. For real sequences $\{f_k\}$ and $\{g_k\}$, we write $f_k\sim g_k$ if $f_k/g_k\to c$ as $k\to\infty$ for some $c>0$. For the Gamma function, we have that $\Gamma(k+\zeta)/\Gamma(k)\sim k^\zeta$, and it follows that $r_i^{\sss(k)}\sim k^{-\tau_i}$ as $k\to\infty$, where
\begin{equation}\label{eq:PA_exp}
\tau_i:=1+\frac{p_i}{p_i+a_i}=2+\frac{p_i}{p_i\theta_i+p_{i^c}(1-\theta_{i^c})}.
\end{equation}
The exponent $\tau_i$ is hence a natural function of the fraction $p_i$ of type $i$ vertices and the fraction $a_i$ of vertices connecting to type $i$ vertices, and can be viewed as quantifying the preferential attachment effect per type $i$ vertex. For the degree $D(t)$ of a randomly chosen vertex, we obtain that
$$
\lim_{t\to\infty}\PP(D(t)=k)=p_1r^{\sss(k)}_1+p_2r^{\sss(k)}_2\sim k^{-\mbox{min}\{\tau_1,\tau_2\}}.
$$
How are the degree distributions of the types affected by the various parameters involved? We first note that, when $p_i=\theta_i=1$, the model reduces to the standard model with only one vertex type. The exponent then becomes 3 which is in agreement with previous results. We obtain the same exponent 3 for both types for any fixed value of $p_i$ when $\theta_i=\theta_{i^c}=1$. The population is then totally homophilic, in the sense that a vertex always connects to its own type, and is thereby divided in two parts, each evolving as the standard model. In a totally heterophilic population, with $\theta_i=\theta_{i^c}=0$, the exponent for type $i$ is $2+\frac{p_i}{p_{i^c}}$. The type that occupies the smallest fraction of the population hence has a degree distribution with heavier tail, as one would expect. In the limit when $\theta_i\to 1$ and $1-\theta_{i^c}=1$, so that both types connects only to type $i$ vertices, we get $\tau_i=2+p_i$ and $\tau_{i^c}=\infty$. Indeed, all type $i^c$ vertices then have degree 1, while type $i$ obeys a power law with heavier tail than in the standard model. If $a_i$ increases while $p_i$ is fixed, then the exponent decreases. Conversely, if $p_i$ increases while $a_i$ is fixed, then the exponent increases.

When $p_1=p_2=1/2$, we get $\tau_i=2+\frac{1}{1+\theta_i-\theta_{i^c}}$. In this case $\tau_1<\tau_2$ if and only if $\theta_2<\theta_1$, which is equivalent to having $a_1>1/2$. Type 1 hence has a heavier tail if and only if the probability that a new vertex attaches to type 1 is larger than for type 2.

\subsection*{Simulations}

We have simulated the model for a few different sets of parameter values and estimated the resulting exponents. The simulated instances of the model are described in Table \ref{tab:PA_cases}, and the code is available online; see \cite{code}. For each set of parameter values, we have generated one graph with $t=10^9$ vertices and base our estimates on this. In this subsection we hence fix $t=10^9$ and omit the dependence on $t$ in the notation. Note that, if $\PP(D=k)\sim k^{-\tau}$, then $\PP(D\geq k)\sim k^{-\gamma}$, where $\gamma=\tau-1$. Throughout we have estimated tail exponents $\gamma$ based on empirical distribution functions rather than exponents $\tau$ based on empirical probability density functions, since this gives rise to more stable estimates.

\begin{figure}
\centering \mbox{\subfigure[Case I.]{\includegraphics[height=3.5cm]{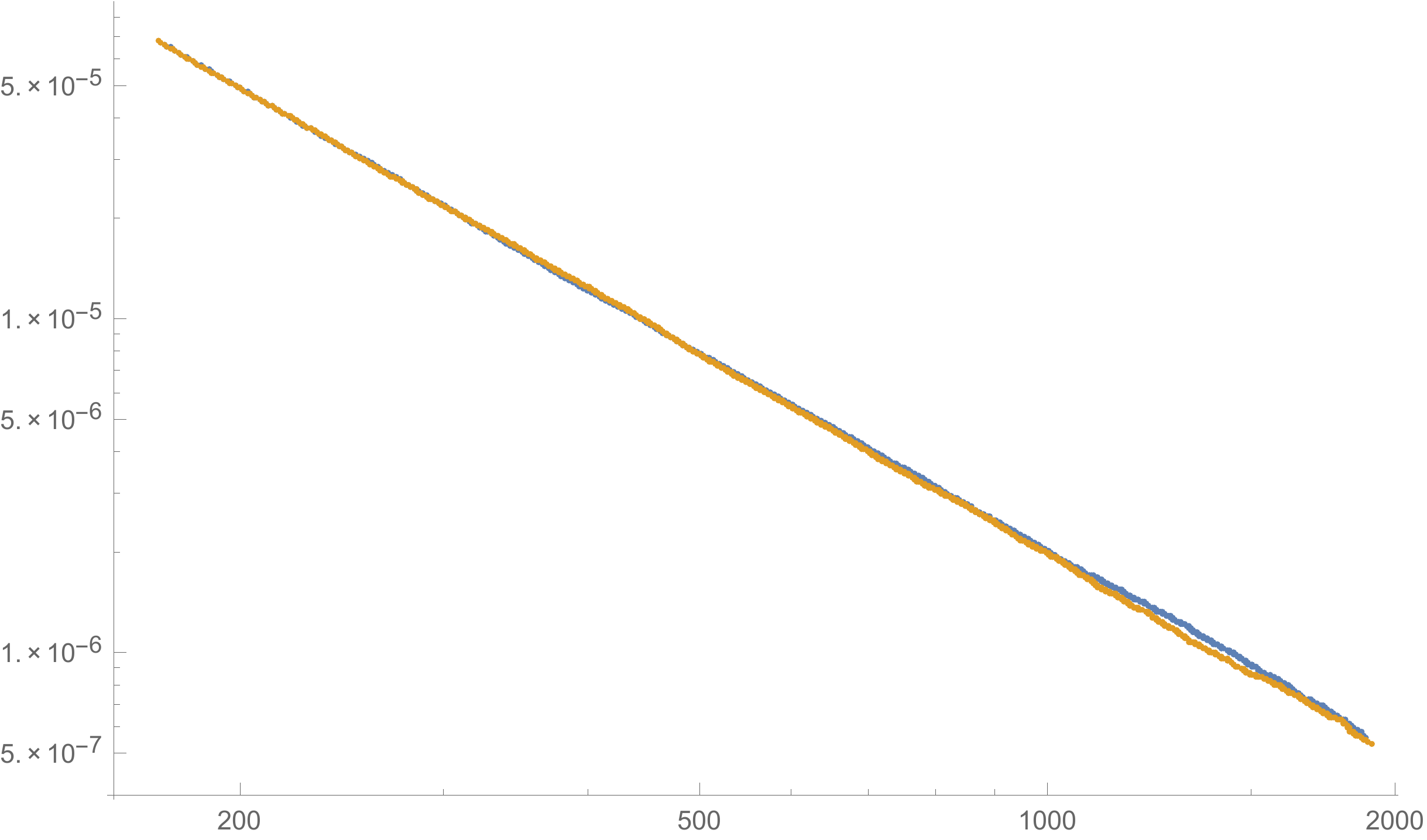}}}
\centering \mbox{\subfigure[Case V.]{\includegraphics[height=3.5cm]{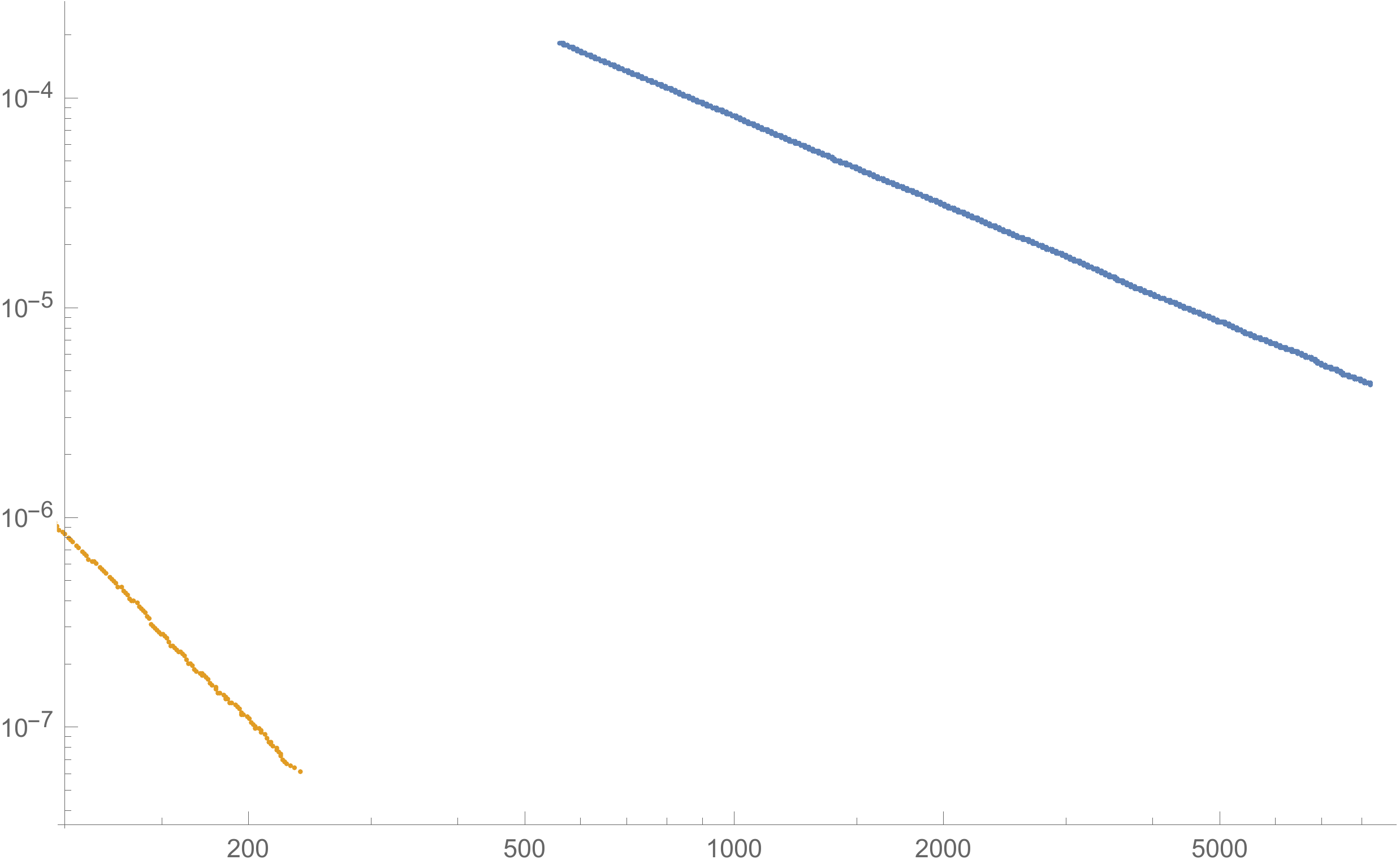}}}
\caption{The tail $N^{\sss(\geq k)}_i$ of the distribution function plotted against $k$ on a log-log scale for type $i=1$ (blue) and type $i=2$ (yellow).}\label{fig:PA_loglog}
\end{figure}

Let $N^{\sss(\geq k)}_i$ denote the number of type $i$ vertices with degree at least $k$, and write $\bar{N}^{\sss(\geq k)}_i=N^{\sss(\geq k)}_i/N_i$. We first note that the simulations strongly suggest that the degree sequences indeed obey power-laws. Figure \ref{fig:PA_loglog} shows log-log plots of $\bar{N}^{\sss(\geq k)}_i$ against $k$ for case I and V, and both pictures reveal straight lines for both types. Plots for the other cases also indicate linear relations. The tail exponent in the power-law is given by the slope of the line. Table \ref{tab:PA_exp} contains estimated values $\hat{\gamma_i}$ of the tail exponents along with the analytical values $\gamma_i=\tau_i-1$, with $\tau_i$ defined in \eqref{eq:PA_exp}. The exponents are estimated using standard methods for linear regression. Throughout, the estimated exponents agree well with the analytical values. For the symmetric cases I and II, the agreement is very good. For the more asymmetric cases III-V, the agreement is not as good. In these cases, the maximal degree of type 2 is only of the order 100, which might not be sufficient to give accurate estimates.

\begin{table}
\begin{center}
\begin{tabular}{|l | l |l|}
\hhline{---} Case & Parameters & Description \\ \hline
\multirow{2}{*}{I} & $p_1=0.5$ & Equal proportions\\
& $\theta_1=\theta_2=0.8$ & Both types prefer its own type\\ \hline
\multirow{2}{*}{II} & $p_1=0.5$ & Equal proportions\\
& $\theta_1=\theta_2=0.2$ & Both types prefer the opposite type\\ \hline
\multirow{2}{*}{III} & $p_1=0.5$ & Equal proportions\\
& $\theta_1=0.8$,$\theta_2=0.2$ & Both types prefer type 1\\ \hline
\multirow{2}{*}{IV} & $p_1=0.1$ & Type 1 minority\\
& $\theta_1=0.8$,$\theta_2=0.2$ & Both types prefer type 1\\ \hline
\multirow{2}{*}{V} & $p_1=0.2$ & Type 1 minority\\
& $\theta_1=\theta_2=0.2$ & Both types prefer the opposite type\\ \hline
\end{tabular}
\end{center}
\caption{Simulated instances of the two-type preferential attachment model}
\label{tab:PA_cases}
\end{table}

\begin{table}
\begin{center}
\begin{tabular}{|c |l| l |l |l |}
\hhline{-----} Case & $\hat{\gamma_1}$ & $\gamma_1$ & $\hat{\gamma_2}$ & $\gamma_2$ \\ \hline
I & 1.9802 & 2 & 1.9895 & 2 \\
II & 1.9716 & 2 & 2.0022 & 2 \\
III & 1.6101 & 1.6250 & 3.4266 & 3.5 \\
IV & 1.3229 & 1.125 & 5.0564 & 5.5 \\
V & 1.3990 & 1.2941 & 3.4067 & 3.5 \\ \hline
\end{tabular}
\end{center}
\caption{Estimated and analytical tail exponents.}
\label{tab:PA_exp}
\end{table}

As for the values of the exponents, we note that, in Case I and II, the tail exponent is the same for type 1 and type 2 and equals the exponent in the standard one-type preferential attachment model. This is natural in view of the symmetry of the population composition and the preferences. In Case III, when both types prefer type 1, the type 1 exponent becomes smaller, while the type 2 exponent becomes larger, as expected. Case IV is similar, but here the preferred type 1 is a minority. The types attract the same fractions of the new arrivals as in Case III, but type 1 (2) occupies a smaller (larger) fraction of the population, meaning that the type 1 (2) exponent becomes smaller (larger). Finally, in Case V, both types prefer the opposite type and type 1 is a minority. Type 1 then attracts a larger fraction of the new arrivals and thereby gets a smaller exponent than type 2.

\begin{figure}
\centering \mbox{\subfigure[Case I.]{\includegraphics[height=6cm]{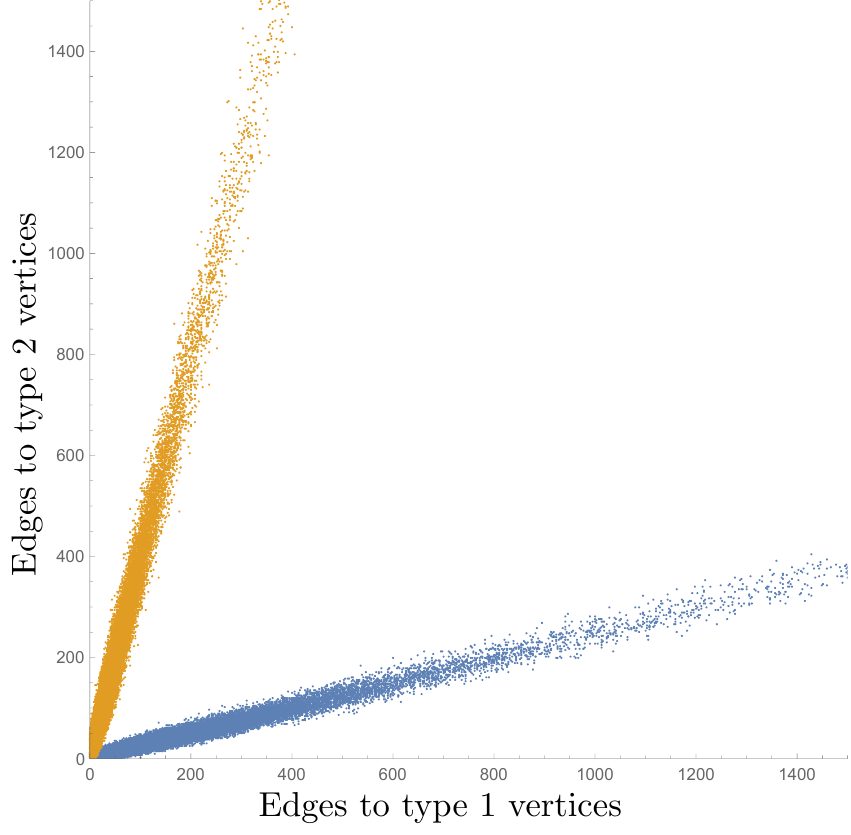}}}\hspace{1cm}
\centering \mbox{\subfigure[Case II.]{\includegraphics[height=6cm]{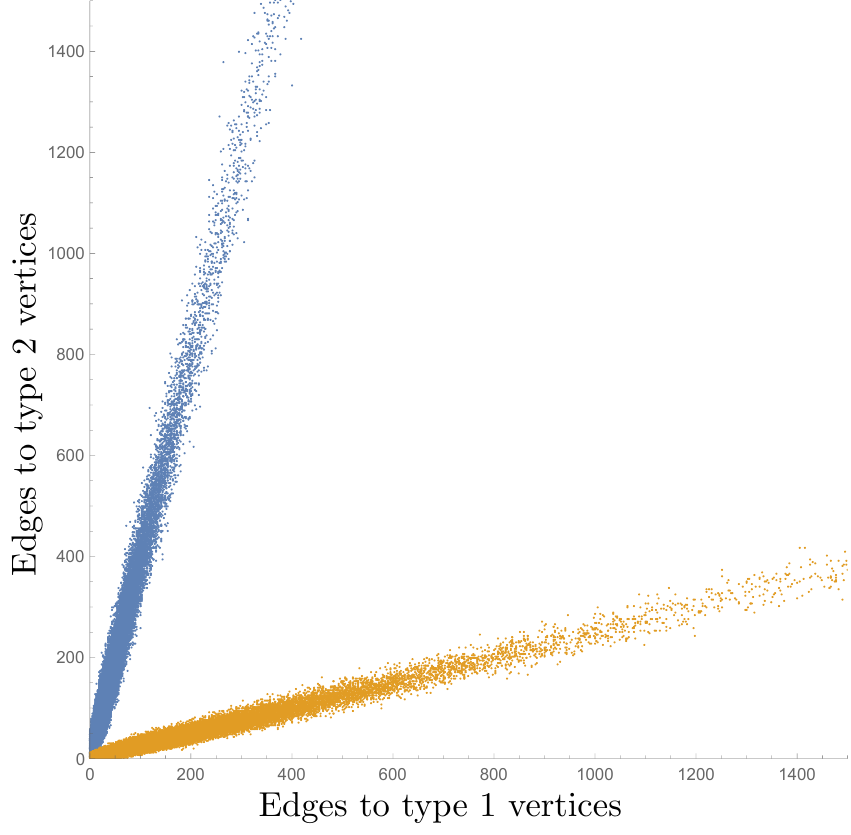}}}\par
\centering \mbox{\subfigure[Case V.]{\includegraphics[height=5cm]{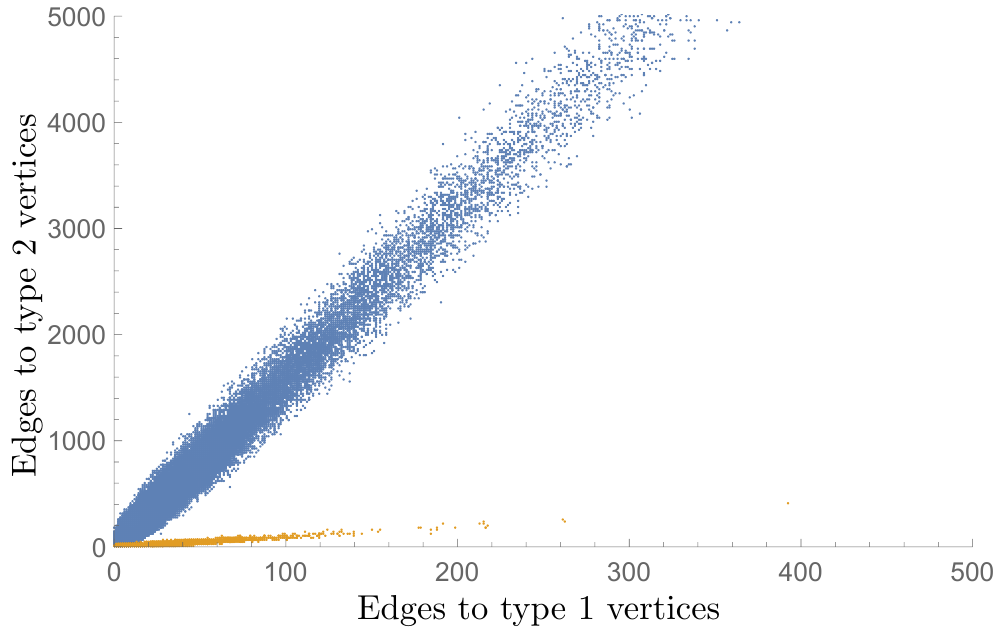}}}
\caption{Scatter plots of type 1 and type 2 degree (that is, number of edges connecting to type 1 and type vertices, respectively) for type 1 vertices (blue) and type 2 vertices (yellow).}\label{fig:PA_corr}
\end{figure}

\begin{table}
\begin{center}
\begin{tabular}{|c|  c|c |c | c|}
  \hline
  {Case}  & $\bar N_{1\to 1}$ & $\bar N_{1\to 2}$  & $\bar N_{2\to 1}$ & $\bar N_{2\to 2}$ \\
  \hline
  I & 1.5968 & 0.3991 & 0.3993 & 1.5972 \\
  II & 0.3993 & 1.597 & 1.5969 & 0.3992 \\
  III & 1.5689 & 0.9688 & 1 & 0.4 \\
  IV & 1.2017 & 3.8155 & 0.8222 & 0.4 \\
  V & 0.3534 & 3.2562 & 1 & 0.4 \\
    \hline
\end{tabular}
\end{center}
\caption{Average degree of the vertices split per type.}
\label{tab:PA_av_per_type}
\end{table}

We have also investigated the degrees split per type. Write $D_{s\to i}$ for the number of type $i$ vertices that vertex $v_s$ is connected to. We will refer to this as the type $i$ degree of vertex $v_s$. We first note that, for each $s$, the pair $\{D_{s\to 1}, D_{s\to 2}\}$ should be positively correlated, since the preferential attachment rule is based on the total degree and a large type $i$ degree hence helps in attracting new vertices also of the other type. That this is indeed the case is illustrated in Figure \ref{fig:PA_corr}, with scatter plots of the type 1 and the type 2 degrees for both types for Case I, II and V. The plots reveal that vertices with a large type 1 degree throughout also tend to have a large type 2 degree. For Case I, we also see that both vertex types tend to have a much larger degree to its own type than to the opposite type while, for Case II, it is the other way around. For Case V, type 1 vertices tend to have a large type 2 degree, while type 2 vertices have much smaller degrees overall (note the different scales on the axes).

Let $\bar{N}_{i\to j}$ denote the average type $j$ degree for the type $i$ vertices. The total expected degree of type $i$ vertices to type $i_c$ is given by $t[p_i(1-\theta_i)+p_{i^c}(1-\theta_{i^c})]$ and hence $\bar{N}_{i\to i^c}$ should be close to $[p_i(1-\theta_i)+p_{i^c}(1-\theta_{i^c})]/p_i$. Similarly, the total expected degree of type $i$ connecting to type $i$ is $2tp_i\theta_i$, so $\bar{N}_{i\to i}$ should be close to $2\theta_i$. Table \ref{tab:PA_av_per_type} contains values of $\bar{N}_{i\to j}$ from the simulations, and it can easily be checked that the values are very close to the predicted ones.

Finally, we look at the tail exponents for the degrees split per type. Write $N^{\sss(\geq k)}_{i\to j}$ for the number of type $i$ vertices that are connected to at least $k$ type $j$ vertices. It is reasonable to expect that also $N^{\sss(\geq k)}_{i\to j}$ obey power-laws. The simulations indicate that this is indeed that case and Table \ref{tab:PA_exp_per_type} shows the estimated associated tail exponents $\hat{\gamma}_{i\to j}$. We note that, for each type, the two exponents -- corresponding to edges leading to the same and the opposite type, respectively -- seem to coincide and be identical to the exponent of the total degree.  One should indeed expect that these exponents are the same, since the preferential attachment mechanism is based on total degree: a given vertex of type $i$ attracts an arriving vertex that has decided to connect to type $i$ according to the same quantity (total degree) regardless of whether the new arrival is of the same or the opposite type.

\begin{table}
\begin{center}
\begin{tabular}{|c |l| l |l |l |}
\hhline{-----} Case & $\hat{\gamma}_{1\to 1}$ & $\hat{\gamma}_{1\to 2}$ & $\hat{\gamma}_{2\to 2}$ & $\hat{\gamma}_{2\to 1}$ \\ \hline
I & 1.9873 & 2.0410 & 1.9926 & 2.0563 \\
II & 2.0154 & 1.9741 & 1.9637 & 2.0062 \\
III & 1.6154 & 1.6141 & 3.5387 & 3.5556 \\
IV & 1.3535 & 1.3230 & 4.9132 & NA \\
V & 1.3827 & 1.3987 & 3.4713 & 3.5191 \\ \hline
\end{tabular}
\end{center}
\caption{Estimated tail exponents for degrees split per type.}
\label{tab:PA_exp_per_type}
\end{table}

\section{Discussion}

We have investigated the behavior of three standard network models when the population is divided in two groups with potentially different connection propensity between and within groups. For two-type versions of the Erd\H{o}s-R\'{e}nyi graph and the configuration model we have seen that the critical parameter for the occurrence of a giant component decreases as the vertices become more prone to connect to vertices of opposite type. The size of the giant component may not be described by the critical parameter, in the sense that we were able to find parameter values such that the component size in our simulations increases while the critical parameter decreases. Throughout, the component size has been studied via simulations, but it may well be possible to obtain analytical results. The two-type Erd\H{o}s-R\'{e}nyi graph for instance is a special case of the model in \cite{BJR} and, in \cite{BJR}, there are precise results on many aspects of the model, including equations for the component size. We have seen that combining a subcritical and a supercritical type leads to an increasing component size if the supercritical type is strong enough and a decreasing component size in other cases. Can it be quantify how strong the supercritical type has to be in order for the component size to be increasing? Also, we have studied how properties of the model change as the graph goes from a homophilic connection pattern to a heterophilic one, but the parameters of the model can of course be tuned in many other ways; see Example 2.5 for a brief example.

For the two-type configuration model, a rigorous analysis of the phase transition could presumably be carried out along the same lines as in \cite{Ball_Sirl}, where a very similar model is treated. An interesting feature to study then is distances between vertices. For the standard configuration model, it is well known that distances are of order $\log n$ if the degree distribution has finite variance and of order $\log\log n$ if the variance is infinite. What about distances in the two-type model if one of the types has a degree distribution with finite variance and the other a degree distribution with infinite variance? Does the answer depend on the proportion of types and how the edges are allocated between/within the types?

For a two-type version of the standard preferential attachment model, we have derived expressions for the degree exponents of the respective types and confirmed these by aid of simulations. In our version of the model, new vertices first decide which type to connect to and then chooses a vertex of the selected type preferentially according to total degree. There are of course many other ways of formulating a two-type model based on preferential attachment, for instance a new vertex could select its partner according to a function that takes both type and degree into account simultaneously. The degree used for the preferential attachment could also be type-specific, so that a new vertex could give different weights to degree corresponding to neighbors of its own type and opposite type, respectively, when deciding where to connect.

\end{document}